%% file: 2dnear.tex
\documentclass[preprint, 10pt]{elsarticle}

\input{preambles.tex}
\usepackage{lineno}

\newcommand{\mcaption}[2]{\caption{\small \em #1}\label{#2}}

\begin{document}

\title{High-volume fraction simulations of two-dimensional vesicle
  suspensions}

\author[ut]{Bryan Quaife} \ead{quaife@ices.utexas.edu}
\author[ut]{George Biros}\ead{gbiros@acm.org}
\address[ut]{Institute of Computational Engineering and Sciences,\\
  The University of Texas at Austin, Austin, TX, 78712.}

\begin{abstract} 
We consider numerical algorithms for the simulation of the rheology of
two-dimensional vesicles suspended in a viscous Stokesian fluid.  The
vesicle evolution dynamics is governed by hydrodynamic and elastic
forces. The elastic forces are due to local inextensibility of the
vesicle membrane and resistance to bending.  Numerically resolving
vesicle flows poses several challenges.  For example, we need to
resolve moving interfaces, address stiffness due to bending, enforce
the inextensibility constraint, and efficiently compute the
(non-negligible) long-range hydrodynamic interactions.

Our method is based on the work of {\em Rahimian, Veerapaneni, and
Biros, ``Dynamic simulation of locally inextensible vesicles suspended
in an arbitrary two-dimensional domain, a boundary integral method'',
Journal of Computational Physics, 229 (18), 2010}.  It is a boundary
integral formulation of the Stokes equations coupled to the interface
mass continuity and force balance.  We extend the algorithms presented
in that paper to increase the robustness of the method and enable
simulations with concentrated suspensions.

In particular, we propose a scheme in which both intra-vesicle and
inter-vesicle interactions are treated semi-implicitly.  In addition we
use special integration for near-singular integrals and we introduce a
spectrally accurate collision detection scheme.  We test the proposed
methodologies on both unconfined and confined flows for vesicles whose
internal fluid may have a viscosity contrast with the bulk medium.  Our
experiments demonstrate the importance of treating both intra-vesicle
and inter-vesicle interactions accurately.
\end{abstract}


\maketitle

\section{Introduction\label{s:intro}}
\input introduction.tex

\section{Formulation\label{s:formulation}} 
\input formulation.tex

\section{Method\label{s:method}} 
\input method.tex

\section{Computing Local Averages of Pressure and Stress\label{s:aver}}
\input aver.tex

\section{Results\label{s:results}} 
\input results.tex

\section{Conclusions\label{s:conclusions}}
\input conclusions.tex

\begin{appendices}
\section{Error estimates for near-singular integration \label{A:AppendixA}} 
\input appen1.tex

\section{Jumps in pressure and stress \label{A:AppendixB}} 
\input appen2.tex

\section{Variable curvature formulation \label{A:AppendixC}}
\input appen3.tex

\end{appendices}

\bibliographystyle{plainnat} 
\bibliography{refs}
\biboptions{sort&compress}
\end{document}

%% file: preambles.tex
\usepackage{palatino}
\usepackage[usenames]{color}
\usepackage{epsfig}
\usepackage[fleqn,reqno]{amsmath}
\usepackage{amsfonts,amsthm,bm}
\usepackage{amssymb}
\usepackage{mathrsfs}
\usepackage{stmaryrd}
\usepackage{subfigure}
\usepackage[titletoc]{appendix}
\usepackage{tabularx,booktabs}
\usepackage{graphics,caption}
\usepackage{fancybox}
\usepackage[top=1.2in,bottom=1.2in,left=1in, right=1in]{geometry}
\usepackage{array}
\usepackage{multirow}
\usepackage{wrapfig}
\usepackage{algorithmic,algorithm}
\usepackage{paralist}
\usepackage{ifthen}
\usepackage{enumitem}
\usepackage{tikz,pgfplots,filecontents}
\usepackage{mdframed}

\usepackage[pagebackref=false,bookmarks=false]{hyperref} 

\definecolor{dblue}{rgb}{0.03,0.3,0.62}
\definecolor{dorange}{rgb}{1,0.55,0}

\hypersetup{
  bookmarksnumbered=true,
  bookmarksopen=false,
  hypertexnames=false,      
  breaklinks=true,          
  unicode=false,
  pdffitwindow=true,        
  pdfnewwindow=true,        
  colorlinks=true,         
  linkcolor=dblue,
  anchorcolor=red,
  citecolor=dorange,
  filecolor=magenta,
  urlcolor=dblue,
  pdfstartview = FitH,
  pdfkeywords = {},
  pdfcreator = {LaTeX with hyperref package}
}

\definecolor{sblue}{cmyk}{0.98,0.13,0,0.43} 
\definecolor{sblue}{cmyk}{0.98,0.13,0,0.43} 

\def\gap{\hspace*{.2in}}

\newcommand{\bigO}{\mathcal{O}}

\newcommand{\p}{\partial}

\newcommand{\xx}{{\mathbf{x}}}
\newcommand{\yy}{{\mathbf{y}}}
\newcommand{\ff}{{\mathbf{f}}}
\renewcommand{\SS}{{\mathcal{S}}}
\newcommand{\BB}{{\mathcal{B}}}
\newcommand{\DD}{{\mathcal{D}}}
\newcommand{\EE}{{\mathcal{E}}}

\newcommand{\pderiv}[2]{\frac{\partial #1}{\partial #2}}
\newcommand{\cc}{{\mathbf{c}}}
\newcommand{\uu}{{\mathbf{u}}}

\newcommand{\rr}{{\mathbf{r}}}
\newcommand{\nn}{{\mathbf{n}}}
\newcommand{\eeta}{{\boldsymbol\eta}}

\newcommand{\ttau}{{\boldsymbol\tau}}
\newcommand{\ssigma}{{\boldsymbol\sigma}}
\newcommand{\llambda}{{\boldsymbol\lambda}}
\newcommand{\grad}{{\triangledown}}
\DeclareMathOperator*{\argmin}{\arg\!\min}

\newcolumntype{C}{>{\centering\arraybackslash} m{2.5cm}}

%% file: introduction.tex
Vesicles are deformable capsules filled with a viscous fluid. The
hydrodynamics of vesicles that are suspended in a viscous fluid
(henceforth, ``{\em vesicle flows}'') play an important role in many
biological phenomena~\cite{kraus1996,seifert}.  For example, they are
used experimentally to understand properties of
biomembranes~\cite{sackmann1996} and red blood
cells~\cite{noguchi2005,pozrikidis1990,ghigliotti-biros-e10,kaoui-tahiri-biros-misbah-e11,misbah2006}.
Here we discuss numerical algorithms for simulating the motion of
vesicles in Stokesian fluids in two dimensions.  Many features of
algorithms designed for particulate flows of vesicle suspensions find
applications to other deformable particulate flows, for example
deformable capsules, bubbles, drops, and elastic filaments.  Vesicle
flows are characterized by large deformations, the local
inextensibility of a vesicle's membrane, the conservation of enclosed
area due to the incompressibility of the fluid inside the vesicle, and
the stiffness related to tension and bending forces.  Efficient and
accurate numerical methods need to address these issues.  It is known
that rheology can be sensitive to elastic instabilities that need to be
resolved quite accurately in long-horizon
simulations~\cite{kaoui-biros-misbah09}. 

In~\cite{pozrikidis2001a}, Pozrikidis presents integral equation
formulations for different types of particulate and interfacial viscous
flows.  In line with our work, integral equation methods have been
applied to vesicle suspensions and other similar capsules in two- and
three-dimensions by many
groups~\cite{veerapaneni-e11,rah:vee:bir,shravan,fre:ore2011,fre:zha2010,zha:isf:ols:fre2010,zha:sha2011a,zha:sha2013a,zha:sha2013b,zha:sha:nar2012},
to name a few.  We focus on algorithms for vesicle flows from our
group's previous work~\cite{shravan}, where numerical algorithms for a
boundary integral formulation of the vesicle hydrodynamics are
presented.  In this formulation, a Stokes single-layer potential is
coupled with the vesicle membrane forces and results in an
integro-differential equation that is constrained by the local
inextensibility constraint (see Section~\ref{s:formulation} for
details).  In~\cite{rah:vee:bir}, these methods were extended to
confined flows and problems with vesicles whose enclosed fluid has a
viscosity contrast with the bulk fluid.  In~\cite{veerapaneni-e11}, we
developed algorithms for three-dimensional vesicle flows.  To summarize,
the main components of the formulation we have used in our prior work
are the following: (1) We proposed a time-stepping method in which the
inter-vesicle interactions are treated explicitly and the self-vesicle
interactions are treated semi-implicitly to remove the bending and
tension stiffness; (2) We proposed a spectrally accurate discretization
in space for differentiation and integration of smooth functions and
weakly singular integrals; (3) We proposed effective preconditioners for
the self-vesicle interactions; and (4) We used the fast multipole method
to accelerate the integral operators that capture the long range
hydrodynamic interactions.

\paragraph{Contributions}
In this paper, we continue our efforts towards efficient algorithms for
the simulation of vesicle hydrodynamics using boundary integral
formulations.  We focus on specific issues that are common when
simulating high concentration vesicle flows.  In particular, we look
into stiffness due to {\em inter-vesicle} dynamics, spectrally accurate
collision detection, near-singular integration, and calculating
quantities like pressure and stress fields that can be used to define
coarse-grained variables.  To make things somewhat more precise, let
$\xx_{p}$ be the shape parametrization of the $p^{\mathit th}$ vesicle
(e.g., Lagrangian points or Fourier coefficients), and let $\uu_{p}$ be
its velocity. Then the evolution equation for the position of the
boundaries of $M$ vesicles can be summarized by
\begin{align}
 \frac{d \xx_p}{d t} =\uu_\infty(\xx_p) + \sum_{q=1}^M \uu_{p}(\xx_{q}), \quad
p=1,\ldots,M,
\label{e:highlevel}
\end{align}
where $\uu_{p}(\xx_{q})$ is the velocity induced on vesicle $p$ due to
forces on vesicle $q$ (via hydrodynamic interactions) and
$\uu_\infty(\xx_p)$ is an external velocity field that drives the flow.
The term $\uu_{p}(\xx_{p})$ is the self-interaction term for the vesicle
$p$ and $\uu_{p}(\xx_{q})$, $q \neq p$, are the inter-vesicle
interactions.

\begin{itemize}
\item {\bf Implicit inter-vesicle interactions.} Most research groups
working on vesicle simulations discretize~\eqref{e:highlevel} using an
explicit time-stepping method. However,~\eqref{e:highlevel} is quite
stiff due to tension and bending.  For this reason, in~\cite{shravan}
we introduced time-marching schemes that treat the self-interactions
semi-implicitly.  For example, a first-order accurate in time scheme
(with $\uu_\infty=0$) reads
\begin{align*}
  \frac{\xx_p^{t+dt}-\xx_p^t}{dt} = \tilde{\uu}_p(\xx_p^{t+dt}) + 
    \sum_{\substack{q=1 \\ q \neq p}}^{M} \uu_{p}(\xx_{q}^{t}), 
      \quad p=1,\ldots,M. 
\end{align*}
The velocity $\tilde{\uu}_{p}(\xx_{p}^{t+dt})$ indicates a linearized
semi-implicit approximation of the vesicle self-interaction term
$\uu_{p}(\xx_{p})$.  Thus, computing $\xx_{p}^{t+dt}$ costs one linear
solve per vesicle.  This is sufficient to reduced bending and tension
stiffness but can be problematic for concentrated suspensions since
$\uu_{p}(\xx_{q})$ also introduces stiffness when vesicles $p$ and $q$
are near each other.  In this paper, we introduce a semi-implicit
scheme of the form
\begin{align*}
  \frac{\xx_p^{t+dt}-\xx_p^t}{dt} = \sum_{q=1}^M 
    \tilde{\uu}_p(\xx^{t+dt}_{q}), \quad p=1,\ldots,M,
\end{align*}
that treats all coupling terms in a linear semi-implicit way.  Thus,
one needs to solve a coupled system for all vesicles at every time
iteration.  We examine the stability of this scheme for first- and
second-order time-stepping and we compare it with alternative
time-marching methods.

\item {\bf Near-singular integration.} When the vesicles come closer
together, the hydrodynamic interaction terms (e.g., $\uu_{p}(\xx_{q})$)
require the evaluation of near-singular integrals.  This can be quite
expensive computationally.  Here we use the two-dimensional analogue of
the scheme introduced in~\cite{ying-biros-zorin06} which allows a
fifth-order accurate evaluation of near interactions using optimal work.

\item {\bf Spectral collision detection.} Another issue related to
concentrated suspensions is collision detection.  While it is
physically impossible for vesicles to collide with one another or with
solid walls in finite time, it is important to detect collisions caused
by numerical errors.  There exist many efficient algorithms for
computing intersections for polygonal domains~\cite{jimenez-e13}.
Since our spatial discretization is done in high-order accuracy and the
vesicle boundaries are $C^\infty$ curves that are spectrally resolved,
we propose a collision detection scheme that is based on potential
theory and it is ideally suited for smooth geometries.  The detection
can be done in linear time in the number of degrees of freedom used to
represent the vesicle boundaries and requires one evaluation of the
free-space Laplace potential.  The method is described in
Section~\ref{s:collision}.

\end{itemize}
We test these methodologies in a variety of flows and we study the
stability and accuracy of the overall algorithm.  As a secondary
contribution, we also provide a scheme for computing average pressures
and stresses in regions of interest, described in
Section~\ref{s:aver}.  The results of our tests are described in
Section~\ref{s:results} in which we test problems with viscosity
contrast, confined flows, and concentrated suspension flows.

\paragraph{Limitations} 
The main limitation is that the method is developed in two dimensions.
Although several flows can be described to good accuracy under a
two-dimensional approximation, a lot of interesting phenomena occur
only in three-dimensional, especially for concentrated suspensions.
The extensions that we present in this paper do not rely on
two-dimensions, and using these methods in three-dimensions is a
problem of implementation.  In particular, the spectral collision
detection scheme only relies on a standard potential theory result that
holds in three-dimensions, and the fast multiple method.  The time
integrators also naturally extend to three-dimensions, but more work is
required to precondition the linear operators that appear upon time
discretization.  The near-singular integration details have been
discussed for generic surfaces in~\cite{ying-biros-zorin06}, but
several optimizations are possible for shapes represented by spherical
harmonics, which we use in our three-dimensional schemes.

Another major limitation is that we do not use time and space
adaptivity.  Those two are essential for efficient robust solvers that
can be used by non-experts.  Currently, we select the time step size by
a trial and error process.  We are currently working on devising time
and space adaptive methods.

Finally, the method is most suitable for Stokesian bulk fluids.  If
inertial or viscoelastic effects are important, a boundary integral
formulation cannot be used.

\paragraph{Related work} 
There is extensive work on vesicle simulations.  Our list here is by no
means exhaustive, but it includes some of the work that is most relative
to ours.  Pozrikidis offers an excellent review of numerical methods for
interfacial dynamics in a Stokes flow~\cite{pozrikidis2001a}.  Capsules
similar to ours have been simulated in concentrated suspensions
in~\cite{li:poz2002,zha:sha:nar2012,zha:sha2011b,zha:isf:ols:fre2010,fre:zha2010,fre:ore2011}.
Also, in~\cite{pra:riv:gra2012,kum:riv:gra2014,kum:gra2011}, high
concentration suspensions of a different kind of capsule are considered,
but these capsules do not resist bending; therefore, the governing
equations are far less stiff.  While some of these methods use spectral
methods, most of them use an explicit time stepping method which results
in a strict restriction on the time step size.  In order to allow for
larger time steps, a Jacobian-free Newton method outlined
in~\cite{dim:hig1997} was applied to an implicit discretization of
droplets submerged in a Stokes flow~\cite{dim2007}.  While this allows
for larger time steps, it requires the solution of nonlinear equations
which can be computationally expensive, and this technique has not been
tested on vesicle suspensions.  Work that uses a combination of implicit
and explicit methods to study the dynamics of a single vesicle
include~\cite{zha:sha2013a,zha:sha2009}, but both these methods require
multiple solves per time step since they use a predictor-corrector
scheme.  In this paper, we remove the stiffness with a semi-implicit
method which only requires solving one linear equation per time step.
These methods have been applied to a single vesicle
in~\cite{zha:sha2009,sal:mik2012} and are extended to multiple vesicles,
where the only the self-interactions are treated semi-implicitly
in~\cite{rah:vee:bir,shravan,zha:sha2013b}.  However, with multiple
vesicles, we are unaware of any work that couples all the vesicles and
solid walls implicitly and  achieves more than first-order accuracy in
time.

Evaluating layer potentials close to their sources is an active area of
research.  Popular methods to evaluate near-singular integrals include
upsampling or high-order quadrature
rules~\cite{helsing-ojala08,kro1999}, and singularity subtraction or
partitions of unity~\cite{poz1999,fre:zha2010,zha:isf:ols:fre2010}.
Some of these methods claim to achieve up to third-order accuracy, but
they do not report timings or accuracies.  Therefore, the cost required
to achieve the accuracies we desire is unclear.  Moreover, these methods
often depend on the dimension of the problem and nature of the
singularity.  A more recent technique is QBX~\cite{klo:bar:gre:one2012}
which can deliver arbitrary accuracy, but is difficult to implement and
we anticipate it will be too expensive for problems with moving
boundaries.  The method we use naturally extends to three dimensions and
we have successfully used it for a variety of near-singular integrals.

\paragraph{Outline of the paper} 
In Section~\ref{s:formulation}, we summarize the formulation of our
problem. In Section~\ref{s:method}, we discuss the spatio-temporal
discretization including the new scheme that semi-implicitly couples the
vesicle updates, the near-singular evaluation
(Section~\ref{s:near-singular}), and the collision detection
(Section~\ref{s:collision}) .  The average pressure and stress
calculations are in Section~\ref{s:aver} and the results are described
in Section~\ref{s:results}.

%% file: formulation.tex
Let us first define the main variables used to model vesicle flows.
Neglecting inertial forces, the dynamics of vesicle flows is fully
characterized by the position of the interface $\xx(s,t) \in \gamma$,
where $s$ is the arclength, $t$ is time, and $\gamma$ is the membrane of
the vesicle. The position is determined by solving a moving interface
problem that models the mechanical interactions between the viscous
incompressible fluid in the exterior and interior of the vesicle (with
viscosity $\mu$) and the vesicle membrane.  Given $\xx(s,t)$, derived
variables include the fluid velocity $\uu$, the fluid stress $T$, the
pressure $p$, and the membrane tension $\sigma$.  In addition, if $\nn$
is the outward normal to the vesicle membrane $\gamma$, then the stress
jump $\ff(\xx) = \llbracket T \rrbracket \mathbf{n}$, which is a
nonlinear function of the position $\xx$, is equal to the sum of a force
due to the vesicle membrane bending modulus $\kappa_b$ and a force due
to the tension $\sigma$.  For wall-confined flows, additional parameters
are the prescribed wall velocity $\mathbf{U}(\xx,t), \xx \in \Gamma$.

Given these definitions, the equations for a single vesicle flow are
given
\begin{equation}
\label{e:vesicles-pde}
\begin{split}
\mu \grad \cdot (\grad \uu + \grad \uu^{T}) - 
\grad p(\xx) = 0, &\hspace{20pt} \xx \in \Omega\backslash\gamma, \gap
&&\mbox{conservation of momentum},\\
 \grad \cdot \uu(\xx) = 0,  &\hspace{20pt} \xx \in \Omega \backslash
 \gamma, \gap &&\mbox{conservation of mass}, \\
  \xx_{s} \cdot \uu_{s} = 0, &\hspace{20pt} \xx \in \gamma, \gap
  &&\mbox{membrane inextensibility},\\
   \uu(\xx,t) = \dot{\xx}(t), &\hspace{20pt} \xx \in \gamma, \gap   &&\mbox{velocity continuity},\\
\ff(\xx)= -\kappa_{b}\xx_{ssss}
 + (\sigma(\xx) \xx_{s})_{s}, &\hspace{20pt} \xx \in \gamma, \gap  &&\mbox{nonzero stress jump},\\
 \uu(\xx,t) = \mathbf{U}(\xx,t), &\hspace{20pt} \xx \in \Gamma, \gap  &&\mbox{wall velocity}. \\
\end{split}
\end{equation}
The viscosity contrast is taken to be constant inside each vesicle.
However, these values can differ from the viscosity of the exterior
fluid.  We also consider a different stress jump that corresponds to a
prescribed intrinsic curvature for the vesicle membrane.  The
modification to the formulation and results are reported in
Appendix~\ref{A:AppendixC}.  In the case of a problem with $M$ vesicles
with interfaces denoted by $\{\gamma_p\}_{p=1}^M$, we define
\begin{align*}
  \gamma = \mathop{\cup}_{p=1}^M \gamma_p.
\end{align*}
Finally, if $\Omega$ is m-ply connected, we let $\Gamma_{0}$ denote the
connected component of $\Gamma = \Gamma_{0} \cup \Gamma_{1} \cup \cdots
\cup \Gamma_{m}$ that surrounds all the other connected components of
$\Gamma$.

There exist many numerical methods for solving interface evolution
equations like~\eqref{e:vesicles-pde}. Since the viscosity is piecewise
constant with a discontinuity along the interface, we opt for an
integral equation formulation using the Stokes free-space Green's
function. Next, following \cite{rah:vee:bir}, we introduce integral and
differential operators that we will need to
reformulate~\eqref{e:vesicles-pde}.

\subsection{Integral equation formulation}
We give the formulation for the general case in which we have viscosity
contrast between the interior and exterior fluid ($\mu_0$ is the
viscosity of the exterior fluid, $\mu_p$ is the viscosity of the
interior fluid for vesicle $p$, and $\nu_p = \mu_p/\mu_0$), and solid
boundaries with a prescribed velocity.  First we introduce $\SS_{pq}$
and $\DD_{pq}$, the single- and double-layer potentials for Stokes flow,
where the constant factors are chosen so that our formulation is in
agreement with Pozrikidis~\cite[equation 2.2]{pozrikidis2001b}. The
subscripts denote the potentials induced by a hydrodynamic density on
the membrane of vesicle $q$ and evaluated on the membrane of vesicle
$p$:
\begin{align*}
  \SS_{pq}[\ff](\xx) &:= \frac{1}{4\pi\mu_{0}}\int_{\gamma_{q}}\left(
    -\boldsymbol{I} \log \rho
  + \frac{\rr \otimes \rr}{\rho^{2}} \right)\ff\, ds_{\yy},
  && \xx \in \gamma_{p}, \\
  \DD_{pq}[\uu](\xx) &:= \frac{1-\nu_{q}}{\pi}\int_{\gamma_{q}}
    \frac{\rr \cdot \nn}{\rho^{2}}\frac{\rr \otimes \rr}{\rho^{2}}\uu\,
  ds_{\yy}, && \xx \in \gamma_{p},
\end{align*}
where $\rr = \xx - \yy,$ and $\rho = \|\rr\|_{2}$.      
Also, we define
\[\SS_{p} := \SS_{pp}\quad\mbox{and}\quad \DD_{p} := \DD_{pp},\]
to indicate vesicle self-interactions.  Next, we define
\begin{align*}
  \EE_{pq}[\ff,\uu](\xx) &= \SS_{pq}[\ff](\xx) + \DD_{pq}[\uu](\xx), && \xx \in \gamma_{p}, \\
  \EE_{p}[\ff,\uu](\xx) &= \sum_{q=1}^{M} \EE_{pq}[\ff,\uu](\xx), && \xx \in \gamma_{p}.
\end{align*}
$\BB$ is the completed double-layer operator for confined
Stokes flow with density $\eeta$
\begin{align*}
  \BB[\eeta](\xx) = \DD_{\Gamma}[\eeta](\xx) + 
    \sum_{q=1}^{M}R[\xi_{q}(\eeta),\cc_{q}](\xx) + 
    \sum_{q=1}^{M}S[\llambda_{q}(\eeta),\cc_{q}](\xx), 
    \quad \xx \in \gamma \cup \Gamma.
\end{align*}
If $\xx \in \Gamma_{0}$, the rank one modification
$\mathcal{N}_{0}[\eeta](\xx) = \int_{\Gamma_{0}}(\nn(\xx) \otimes
\nn(\yy)) \eeta(\yy)ds_{\yy}$ is added to $\BB$ which removes the
one-dimensional null space of the corresponding integral
equation~\cite{pozrikidis1992}.  The Stokeslets and rotlets are defined
as
\begin{align*}
  S[\llambda_{q}(\eeta),\cc_{q}](\xx) = \frac{1}{4\pi\mu_{0}}
    \left(\log \rho + \frac{\rr \otimes \rr}{\rho^{2}}\right)
    \llambda_{q}(\eeta) 
  \qquad \text{and} \qquad
  R[\xi_{q}(\eeta),\cc_{q}](\xx) = \frac{\xi_{q}(\eeta)}{\mu_{0}}
    \frac{\rr^{\perp}}{\rho^{2}},
\end{align*}
where $\cc_{q}$ is a point inside $\omega_{q}$, $\omega_{q}$ is the
interior of vesicle $q$, $\rr=\xx - \cc_{q}$, and $\rr^{\perp} =
(r_{2},-r_{1})$.  The operator $\BB$ satisfies the jump condition
\begin{align*}
  \lim_{\substack{\xx \rightarrow \xx_{0} \\ \xx \in \Omega}}
    \BB[\eeta](\xx) = -\frac{1}{2}\eeta(\xx_{0}) + 
    \BB[\eeta](\xx_{0}), \qquad \xx_{0} \in \Gamma,
\end{align*}
and the size of the Stokeslets and rotlets are
\begin{align*}
  \llambda_{q,i} = \frac{1}{2\pi} \int_{\gamma_{q}} 
    \eeta_{i}(\yy) ds_{\yy}, \quad i=1,2
  \qquad \text{and} \qquad
  \xi_{q} = \frac{1}{2\pi}\int_{\gamma_{q}} \yy^{\perp}
    \eeta(\yy)ds_{\yy}.
\end{align*}
The inextensibility constraint is written in operator form as
\begin{align*}
  \mathcal{P}[\uu](\xx) = \xx_{s} \cdot \uu_{s}.
\end{align*}

Putting everything together, the integral formulation equation of~\eqref{e:vesicles-pde} is
given by~\cite{rah:vee:bir}
\begin{subequations}
\label{e:ves:dyn}
\begin{align}
  &(1+\nu_{p}) \uu(\xx) = \EE_{p}[\ff,\uu](\xx) + 
  \BB_{p}[\eeta](\xx), \quad &&\xx \in \gamma_{p},
  \label{e:vesicle:dynamic} \\
  &\mathcal{P}[\uu](\xx) = 0, &&\xx \in \gamma_{p},
  \label{e:inextensibility} \\
  &\mathbf{U}(\xx) = -\frac{1}{2}\eeta(\xx) + \EE_{\Gamma}[\ff,\uu](\xx) + 
    \BB[\eeta](\xx), &&\xx \in \Gamma,
  \label{e:BIE}
\end{align}
\end{subequations}
where $\ff$ is the hydrodynamic density, which in our case is the stress
jump on the vesicle interface, which also is referred to as the traction
jump, and is given by 
\begin{align*}
  \ff = -\kappa_{b}\xx_{ssss} + (\sigma \xx_{s})_{s}.
\end{align*}
Since $\uu = d\xx/dt$ and $\ff$ depends on $\sigma$ and
$\xx$, \eqref{e:ves:dyn} is a system of integro-differential-algebraic equations
for $\xx,\sigma$, and $\eeta$.

%% file: method.tex
Following~\cite{rah:vee:bir}, we use a Lagrangian formulation in which
we track the position of material points on $\gamma$ using multistep
implicit-explicit (IMEX) methods~\cite{ascher1995} to
discretize~\eqref{e:ves:dyn} in time. We use Nystr\"{o}m-type
quadrature rules to discretize integral operators and Fourier
differentiation to compute derivatives with respect to the membrane
parametrization.

\subsection{Spatial Discretization\label{s:spaceDisc}} 
Let $\xx(\alpha),$ with $\alpha \in (0,2\pi]$, be a parametrization of
the interface $\gamma_p$, $\{\alpha_k = 2k\pi/n\}_{k=1}^n$ be $n$
material points, and
\begin{align*}
  \xx(\alpha) = \sum_{k = -n/2+1}^{n/2} \hat{\xx}(k) e^{ik\alpha}.
\end{align*}
Using the FFT to compute $\hat{\xx}$, we compute derivatives of $\xx$
with spectral accuracy, assuming that $\gamma$ belongs to $C^{\infty}$.
In particular, we can compute the arclength derivative with spectral
accuracy since 
\begin{align*}
  \pderiv{}{s} = \pderiv{\alpha}{s}\pderiv{}{\alpha} = 
    \frac{1}{\|\partial\xx/\partial\alpha\|}\pderiv{}{\alpha},
\end{align*}
where $s$ is the arclength parametrization.

Since the single-layer potential $\SS$ has a logarithmic singularity, we
use the hybrid Gauss-trapezoid quadrature rule given in Table 8 of
\cite{alpert1999}.  The error of this quadrature rule is $\bigO(h^{8}
\log h)$ for integrands with logarithmic singularities. Letting $\xx_k =
\xx(\alpha_k)$, the quadrature rule is
\begin{align*}
  \mathcal S[\ff](\xx) \approx \sum_{k=1}^{n+m} w_{k}
      S(\xx,\xx_{k}) \ff(\xx_{k})\| \xx_{\alpha,k}\|,
\end{align*}
where $n$ is the number of nodes, $m$ is a number of additional
quadrature nodes, $w_{k}$ are the quadrature weights, $\xx_{k}$ are the
quadrature abscissae, and $\xx_{\alpha}$ is $\partial \xx/\partial
\alpha$.  The first $n$ abscissae are the usual equispaced quadrature
points.  Then, the other $m$, which is determined by the desired order
of convergence for the integral, are additional quadrature abscissas
clustered around the singularity.

The double-layer potential has no singularity in two dimensions since
\begin{align*}
  \lim_{\substack{\xx' \to \xx \\ \xx' \in \gamma}} D(\xx',\xx) = 
    \frac{\kappa(\xx)}{2\pi}({\bf t}(\xx) \otimes {\bf t}(\xx)),
  \quad \xx \in \gamma,
\end{align*}
where $\kappa$ is the curvature and ${\bf t}$ is the unit tangent.
Thus, a composite trapezoid rule will give spectral accuracy since the
integrand is periodic and smooth.

\subsection{Time discretization}
Time derivatives are discretized as
\begin{align*}
  \frac{d\xx}{dt} \approx \frac{\beta\xx^{N+1} - \xx^{0}}{\Delta t},
\end{align*}
where $\xx^{0}$ is a linear combination of previous time steps.
Operators and terms that are treated explicitly are discretized at
$\xx^{e}$ which is also a linear combination of previous time steps.
The simplest IMEX method is IMEX Euler which is given by $\beta=1$,
$\xx^{0} = \xx^{N}$, and $\xx^{e} = \xx^{N}.$ The second-order time
integrator we use is given by $\beta=3/2,$
\begin{align*}
\xx^{0} = 2\xx^{N} - \frac{1}{2}\xx^{N-1}, \hspace{20pt}
\xx^{e} = 2\xx^{N} - \xx^{N-1}.
\end{align*}
To avoid significant time stepping constraints, we treat the bending
term $\ff_{b}$ and tension term $\ff_{\sigma}$ semi-implicitly.  An
approximation for the position and tension of vesicle $p$ at time $N+1$
satisfies
\begin{subequations}
  \begin{align}
    \frac{\alpha_{p}}{\Delta t}(\beta \xx_{p}^{N+1} - \xx_{p}^{0}) &= 
      \SS_{p}^{e} \ff_{p}^{N+1} + 
      \DD_{p}^{e}\uu_{p}^{N+1} + \BB_{p}[\eeta^{N}] 
      + \sum_{\substack{q=1 \\ q \neq p}}^{M} \EE_{pq}^{e}[\ff_{q}^{N},\uu_{q}^{N}], 
      &&\xx \in \gamma_{p},
      \label{e:vesicle:dynamic:exp} \\
      \beta \mathcal{P}^{e}\xx_{p}^{N+1} &= \mathcal{P}^{e}\xx_{p}^{0},
      &&\xx \in \gamma_{p},
      \label{e:inextensibility:exp} \\
      U(\xx) &= -\frac{1}{2}\eeta^{N}(\xx) + \EE_{\Gamma}^{e}[\ff^{N},\uu^{N}](\xx) + 
      \BB[\eeta^{N}](\xx) + \mathcal{N}_{0}[\eeta^{N}](\xx), &&\xx \in \Gamma, 
      \label{e:BIE:exp} \\
      \uu_{p}^{N+1} &= \frac{\beta \xx_{p}^{N+1}-\xx_{p}^{0}}{\Delta t},
      &&\xx \in \gamma_{p},
      \label{e:velocity:exp}
  \end{align}
  \label{e:explicit:method}
\end{subequations}
where $\alpha_{p} = (1+\nu_{p})/2$, and operators with a superscript $e$
are discretized at $\xx^{e}$.  Note that~\eqref{e:explicit:method} has
explicit vesicle-vesicle and vesicle-boundary interactions.  That is, we
are approximating the solution of~\eqref{e:ves:dyn} with the solution of
a problem where the vesicles are fully decoupled from one another and
from the boundary.  To solve~\eqref{e:explicit:method}, we first find
$\eeta^{N}$ by solving~\eqref{e:BIE:exp} and then
solve~\eqref{e:vesicle:dynamic:exp} and~\eqref{e:inextensibility:exp}
for the new position and tension of each vesicle independently of all
the others.

In high-concentration flows, a new source of stiffness arises.  As two
vesicles approach one another, $\EE_{pq}$ becomes increasingly large.
A similar result holds for vesicle-boundary interactions where
$\BB_{p}$ introduces stiffness.  We propose the new method with
semi-implicit vesicle-vesicle and vesicle-boundary interactions (the
differences with \eqref{e:explicit:method} being highlighted in red)
\begin{subequations}
  \begin{align}
    \frac{\alpha_{p}}{\Delta t}
    (\beta \xx_{p}^{N+1} - \xx_{p}^{0}) &= 
      \SS_{p}^{e} \ff_{p}^{N+1} + \DD_{p}^{e}\uu_{p}^{N+1} + 
      \BB_{p}[\eeta^{N\textcolor{red}{+1}}] + \sum_{\substack{q=1 \\ q \neq p}}^{M} 
      \EE_{pq}^{e}[\ff_{q}^{N\textcolor{red}{+1}},
      \uu_{q}^{N\textcolor{red}{+1}}],
      &&\xx \in \gamma_{p},
      \label{e:vesicle:dynamic:imp} \\
      \beta \mathcal{P}^{e}\xx_{p}^{N+1} &= \mathcal{P}^{e}\xx_{p}^{0},
      &&\xx \in \gamma_{p},
      \label{e:inextensibility:imp} \\
      U(\xx) &= -\frac{1}{2}\eeta^{N\textcolor{red}{+1}}(\xx) + 
      \EE_{\Gamma}^{e}[\ff^{N\textcolor{red}{+1}},
      \uu^{N\textcolor{red}{+1}}](\xx) + 
      \BB[\eeta^{N\textcolor{red}{+1}}](\xx)+\mathcal{N}_{0}[\eeta^{N\textcolor{red}{+1}}](\xx), &&\xx \in \Gamma, 
      \label{e:BIE:imp} \\
      \uu_{p}^{N+1} &= \frac{\beta \xx_{p}^{N+1}-\xx_{p}^{0}}{\Delta t},
      &&\xx \in \gamma_{p}.
      \label{e:velocity:imp}
  \end{align}
  \label{e:implicit:method}
\end{subequations}
Unlike in~\eqref{e:explicit:method}, vesicles
in~\eqref{e:implicit:method} are implicitly coupled since the operator
$\EE_{pq}^{e}$ is applied to $\ff_{q}^{N+1}$ and $\uu_{q}^{N+1}$.
Similarly, the operator $\BB_{p}$ implicitly couples the vesicles to the
solid walls.  Equation~\eqref{e:implicit:method} is still an
approximation of~\eqref{e:ves:dyn}, but its benefit is that it is
experimentally a more stable method (see Sections~\ref{s:taylorGreen}
and~\ref{s:couette}).

Since~\eqref{e:implicit:method} is fully coupled, it is too expensive to
solve without preconditioning.  Roughly speaking, the number of GMRES
iterations depends on the sizes of $\EE_{pp}$ and $\EE_{pq}$ and the
latter depends on the inter-vesicle proximity.  We apply preconditioned
GMRES~\cite{saad} to~\eqref{e:implicit:method} where we use the
block-diagonal preconditioner.  The number of GMRES iterations depends
on the minimum inter-vesicle distance, but this preconditioned linear
system experimentally results in mesh-independence\footnote{We
experimented with multigrid preconditioners~\cite{quaife:biros2013} but
mesh-independence was only achieved for computationally expensive
preconditioners.}.

To test the preconditioner, one step of a Taylor-Green flow is done
with the number of points per vesicle ranging from $16$ to $128$
(Table~\ref{t:GMRES:preco}).  We solve~\eqref{e:implicit:method} to a
tolerance of $10^{-8}$ with and without preconditioning.  The required
GMRES iterations for the block-diagonal solver grows almost linearly
and the results are consistent with those reported
in~\cite{rah:vee:bir}.  The implicit solver without preconditioning
requires an unacceptable number of GMRES iterations, but with
preconditioning, we see that the required number of GMRES steps is
independent of $N$ \footnote{Notice, however, that the number of GMRES
iterations will depend on how close the vesicles are to one another as
this determines the overall spectrum structure of the operator. But for
a given flow, the method is independent of $N$.}.  Using the
block-diagonal preconditioner will not result in a speedup unless we
precompute and store the block diagonal preconditioner.  We perform
this precomputation for all the experiments in Section~\ref{s:results}.

\begin{table}[htps]
\centering
\begin{tabular}{>{\centering}m{0.5cm} >{\centering}m{3cm} >{\centering}m{3cm} >{\centering}m{3cm}}
$N$ & Block- \\ Diagonal & Unpreconditioned \\ Implicit & Preconditioned
\\ Implicit \tabularnewline
\hline
16 & 46 & 315 & 18 \tabularnewline
32 & 79 & 341 & 15 \tabularnewline
64 & 164 & 882 & 15 \tabularnewline
128 & 294 & 2117 & 15 
\end{tabular}
\mcaption{The number of GMRES iterations required to take one time step
of a Taylor-Green flow using semi-implicit inter-vesicle
interactions~\eqref{e:implicit:method}.  $N$ is the number of
discretization points per vesicle, the second column is the number of
GMRES iterations to invert the block-diagonal preconditioner, the third
column is the number of unpreconditioned GMRES iterations, and the
final column is the number of preconditioned GMRES iterations.  Notice
that without preconditioning, the number of iterations quickly grows
while the preconditioned linear system is
mesh-independent.}{t:GMRES:preco}
\end{table}
We have also observed experimentally that when vesicles are close,
treating inter-vesicle interactions explicitly with second-order time
stepping is unconditionally unstable.  We have experimented with time
step sizes of varying orders of magnitude and each experiment
inevitably becomes unstable.  However, scheme~\eqref{e:implicit:method}
with second-order time stepping is stable once $\Delta t$ is small
enough for the simulation under consideration.

\subsection{Fast summation}
For high-concentration flows, the most expensive part of the simulation
is evaluating the integral operators $\SS$ and $\DD$.  The single-layer
potential $\SS$ can be written as
\begin{align*}
  \frac{1}{4\pi\mu_{0}} \int_{\gamma}  \left(
    -\boldsymbol{I}\log \rho  + 
    \frac{\rr \otimes \rr}{\rho^{2}} \right)\ff ds = 
    -\frac{1}{4\pi\mu_{0}} \int_{\gamma} (\log \rho) \ff ds +
    \frac{1}{4\pi\mu_{0}} \int_{\gamma}
    \left(\frac{\rr \cdot \ff}{\rho^{2}}\right)\rr ds \\
  =-\frac{1}{4\pi\mu_{0}} \int_{\gamma}(\log \rho) \ff ds + 
    \frac{1}{4\pi\mu_{0}} \int_{\gamma}   
    \frac{\rr}{\rho^{2}} (\xx \cdot \ff) ds_{\yy} - 
    \frac{1}{4\pi\mu_{0}} \int_{\gamma} 
    \frac{\rr}{\rho^{2}} (\yy \cdot \ff) ds_{\yy},
\end{align*}
which can be computed with three Laplace FMM applications.  We also use
the FMM to evaluate the single-layer potential at arbitrary target
locations by assigning a charge of zero to these locations.
Unfortunately, the double-layer potential can not be decomposed into a
sequence of Laplace single- and double-layer potentials.  Future work
involves using the kernel-independent FMM~\cite{ying-biros-zorin-03} to
accelerate evaluating the double-layer potential.

\subsection{Near-singular integration}\label{s:near-singular}
Another difficulty with high-concentration simulations is the
evaluation of layer potentials at points close to a boundary.  Integral
operators $\DD$ and $\SS$ can be approximated with high accuracy using
the trapezoid rule if the target point is sufficiently far from the
source points.  As the target point approaches the vesicle,
more points are required to achieve the same accuracy, and, the
required number of points for a fixed error grows without
bound (Appendix~\ref{A:AppendixA}).  A near-singular integration
strategy that guarantees a uniform accuracy for all target points is
required.  

We adopt the strategy from~\cite{ying-biros-zorin06} since it is not
depend on the nature of the integrand and it is simple to implement.
Suppose that a vesicle $\gamma$ is discretized with $N$ points and that
the resulting arclength term is of size $h$.  Let $d(\xx,\gamma) =
\inf_{\yy \in \gamma} \|\xx - \yy\|$ be the distance from $\xx$ to
$\gamma$.  We define the far zone of $\gamma$ as $\Omega_{1}=\{\xx
\:|\: d(\xx,\gamma) \geq h\}$, and the near zone of $\gamma$ as
$\Omega_{0}=\{\xx \:|\: d(\xx,\gamma) < h\}.$ The near and far zones
can be constructed efficiently using a quadtree-structure, similar to
the near-field and far-field interactions in fast multipole codes.

For $\xx \in \Omega_{1}$, the trapezoid rule with $N^{3/2}$ points
guarantees that the error for $\SS$ is $\mathcal{O}(h^{M/2-2})$ and for
$\DD$ is $\mathcal{O}(h^{M/2-4})$, where the density function (in our
case, either $\ff$ or $\eeta$) belongs to $C^{M}$.  These error bounds
are proved in Appendix~\ref{A:AppendixA}.  For a target point $\xx$ in
the near zone, we first find the closest boundary point, $\xx_{0}$,
using two applications of Newton's method with the initial guess being
the closest discretization point to $\xx$.  We then place $m$
interpolation points at the points $\xx_{j} = \xx_{0} + j\beta h
(\xx-\xx_{0})/\|\xx-\xx_{0}\|$ for $j=0,\ldots,m$, where $\beta$ is
slightly larger than $1$ to guarantee that all the interpolation points
are in the far zone (left plot of Figure~\ref{f:nearsing:conv}).  The
layer potential is evaluated at $\xx_{0}$ using a local interpolant with
$N_{\mathrm{int}}$ points on $\gamma$ and evaluated at $\xx_{j},$
$j=1,\ldots,m$ using the $N^{3/2}$-point trapezoid rule.  Then, a
one-dimensional Lagrange interpolation assigns a value to the layer
potential at $\xx$.  Note that the proposed scheme can become more
accurate by using interpolation schemes that compute derivatives of the
solution and building more accurate polynomial interpolations.  This,
however, would significantly complicate the implementation.  Pseudocode
for this algorithm is given in Algorithm~\ref{a:algorithm2}.  In
Appendix~\ref{A:AppendixA}, we show that the error of our algorithm
applied to the single-layer potential $\SS$ is 
\begin{align*}
  \bigO\left(h^{\min\left(N_{\mathrm{int}}-1,
      m,\frac{M}{2}-2\right)}\right),
\end{align*}
and applied to the double-layer potential $\DD$ is
\begin{align*}
  \bigO\left(h^{\min\left(N_{\mathrm{int}}-1,
      m,\frac{M}{2}-4\right)}\right).
\end{align*}

\begin{algorithm}
\begin{algorithmic} 
\FORALL[Loop over far zone]{$\xx \in \Omega_{1}$} 
\STATE $\uu = \mathtt{trapezoid}(\gamma,\xx,\ff)$
\COMMENT{Apply the $N^{3/2}$-point trapezoid rule}
\ENDFOR 
\FORALL[Loop over near zone]{$\xx \in \Omega_{0}$}
\STATE $\xx_{0} = \argmin d(\xx,\gamma)$
\COMMENT{Find the point on $\gamma$ closest to $\xx$}
\STATE $\uu_{0} = \mathtt{interpolate}
    (\xx_{0},\uu_{1},\ldots,\uu_{N_{\mathrm{int}}})$
\COMMENT{Interpolate $\uu$ at $\xx_{0}$}
\FOR{$j=1,\ldots,m$}
\STATE $\xx_{j} = \xx_{0} + jh\beta(\xx-\xx_{0})/\|\xx-\xx_{0}\|$
\COMMENT{Define Lagrange interpolation points}
\STATE $\uu_{j} = \mathtt{trapezoid}(\gamma,\xx_{j},\ff)$
\COMMENT{Evaluate layer potential at $\xx_{j}$}
\ENDFOR 
\STATE $\uu = \mathtt{interpolate}(\xx,\uu_{0},\ldots,\uu_{m})$
\COMMENT{Compute $\uu$ using Lagrange interpolation}
\ENDFOR 
\RETURN $\uu$ \COMMENT{Return the layer potential at all target points}
\end{algorithmic}
\mcaption{LayerPotential($\gamma$,  $\ff$, $\uu$, $\Omega_{0}$,
$\Omega_{1}$)} {a:algorithm2}
\end{algorithm}

\begin{figure}[htps]
  \begin{center}
    \begin{minipage}{0.45\textwidth}
    \input{near-sing.tikz}
    \end{minipage}
    \begin{minipage}{0.45\textwidth}
    \input{nearSingErrors.tikz} 
    \end{minipage} 
  \end{center} 
  \mcaption{\emph{Left:} Lagrange interpolation points in $\Omega_{1}$
  and $\gamma$ are used to approximate layer potentials at $\xx \in
  \Omega_{0}$.  The value at $\xx_{0}$ uses an interpolant from
  neighboring points on the boundary.  All other interpolation points
  are in $\Omega_{1}$ so that the $N^{3/2}$-point trapezoid rule is
  sufficiently accurate.  \emph{Right:} Here we report convergence
  rates. For the blue and the black curves, the expected order of
  convergence is $5$. We obtain $4.7$ and $5.0$, respectively. As the
  green curve corresponds to a linear flow, Lagrange interpolation gives
  exact values.  Hence, machine accuracy is reached once $N$ is
  sufficiently large.}{f:nearsing:conv}
\end{figure}

We test our near-singular integration with three examples.  In all the
examples, 32 target points are placed in $\Omega_{0}$ for
$N=16,\ldots,512$.  The boundary is a $3:2$ ellipse, the boundary values
are interpolated with $N_{\mathrm{int}}=6$ points, and the Lagrange
interpolation uses $m=6$ points.  Since $m = N_{\mathrm{int}}$, we are
using one more point to interpolate in the direction coming off of
$\gamma$.  We do this because the interpolant that is coming off of the
vesicle is always evaluated near the first interpolation point.  Thus,
it may suffer from the Runge phenomenon (we have observed this
experimentally).  This is not the case for the local interpolation used
to find $\uu(\xx_{0})$ since $\xx_{0}$ is located near the middle of the
local interpolation points.  For the first two examples, we compute the
density function of $\SS$ for the two solutions of the Stokes equations:
\begin{align*}
  &\uu(\xx)=(x,-y), \\
  &\uu(\xx)=\left(-\log|\xx| + \frac{\xx \otimes \xx}{|\xx|^{2}}\right)
  \left[ \begin{array}{c}
    1 \\ 1 
  \end{array}\right].
\end{align*}
For the third example, we pick an arbitrary density function and compare
the results from near-singular integration with the results from an
over-refined trapezoid rule.  The error for all three examples is
reported in~Figure~\ref{f:nearsing:conv}.  For the extensional flow $\uu
= (x,-y)$, the Lagrange interpolation is exact since the flow is linear
and the error decays faster than the other examples.  For the other
examples, the error is close to fifth-order which agrees with the order
of accuracy of the near-singular integration

If we have $M$ vesicles and $N$ points per vesicle, the major costs of
the near-singular integration algorithm are:
\begin{itemize}
  \item \emph{Upsampling each vesicle from $N$ to $N^{3/2}$ points:}
  Using the FFT, this requires $\bigO(MN^{3/2}\log N)$ operations.
  \item \emph{Applying the $N^{3/2}$-point quadrature rule to
  $\Omega_{1}$:} With the FMM, this requires $\bigO(MN^{3/2} + MN)$
  operations, and without the FMM, this requires $\bigO(M^{2}N^{5/2})$
  operations.  
  \item \emph{Evaluating the layer potential in $\Omega_{0}$:} Assuming
  there are $\bigO(1)$ points in the near zone, this requires
  $\bigO(MN^{3/2})$ operations.
\end{itemize}
The algorithm can be further sped up by introducing the intermediate and
far zones~\cite{ying-biros-zorin06}.  The intermediate zone is
$\Omega_{1}=\{\xx \:|\: d(\xx,\gamma) \in (\sqrt{h},h]\}$ and the far
zone is $\Omega_{2}=\{\xx \:|\: d(\xx,\gamma) \geq \sqrt{h}\}$.  Then,
the same accuracy can be achieved by using the $N^{3/2}$-point trapezoid
rule in $\Omega_{1}$ and the $N$-point trapezoid rule in $\Omega_{2}$.
While this would reduce the constant of the complexity, it does not
change the overall order and complicates the implementation.

\subsection{Collision detection}\label{s:collision}
Even with semi-implicit treatment of inter-vesicle interactions, due
to time discretization errors, the vesicles can come so close that they
collide or cross a solid wall.\footnote{The streamlines of Stokes
equations never intersect.  Therefore, it is physically impossible for
vesicles to collide with one another or with solid walls in finite
time.  Lubrication theory can be used to study the forces between
nearly touching vesicles, and therefore compute estimates of the
distance between the vesicles, but we choose not to discuss these
forces here.}  There has been a lot of work in collision
detection~\cite{jimenez-e13,redon-e04}, but they are not spectrally
accurate or require complex surface-surface intersections.  Most
methods construct piecewise approximations to the boundaries to
simplify the detection of collisions.  Since our vesicles and solid
walls are periodic, we can use some nice properties of layer potentials
to create a spectrally accurate collision detector.  We use the classic
potential theory result~\cite{kellogg}
\begin{align}
  I_{\Gamma}(\xx) = \frac{1}{2\pi}\int_{\Gamma}\pderiv{}{\nn_{\yy}}\log|\xx - \yy| ds_{\yy} = 
  \begin{cases}
    1 & \xx \in \Omega, \\ 
    \frac{1}{2} & \xx \in \Gamma, \\ 
    0 & \xx \in \mathbf{R}^{2} \backslash \overline{\Omega}, 
  \end{cases}
  \label{e:layerpotentials}
\end{align}
where $\Omega \subset \mathbf{R}^{2}$ is bounded with a smooth boundary
$\Gamma$, and $\nn_{\yy}$ is the outward pointing normal.  We
apply~\eqref{e:layerpotentials} to the configuration of the vesicles
and the solid walls by computing $I_{X}(\xx) = I_{\gamma_{1}}(\xx) +
\cdots + I_{\gamma_{M}} + I_{\Gamma}$ for all $\xx \in \gamma$.  The
maximum value of $I_{X}(\xx)$ is $3/2$ for confined flows and $1/2$ for
unconfined flows if vesicles have not crossed and have not left the
domain.  The maximum value of $I_{X}(\xx)$ will be at least $3/2$ for
unconfined flows, and at least $5/2$ for confined flows if vesicles
have crossed or left the domain.  We note that $I_{X}(\xx)$ can be
evaluated accurately using our near-singular integration scheme, and in
linear time using the FMM. 

We can use $I_{X}(\xx)$ to create a simple adaptive time-stepping
method.  At each time step, we check for collisions.  If a collision has
occurred, we initialize the simulation at a previous time step with a
time step half the size.  We are currently implementing a more rigorous
adaptive time stepping method that uses ideas from spectral deferred
correction methods~\cite{dut:gre:rok,hua:jia:min,min}.  At each time
step, we estimate the local truncation error and then try to commit a
constant amount of error per time step by adjusting the time step size.
We leave the details and the implementation as future work.  Future work
also includes adaptivity in space.  If vesicles are sufficiently far
from one another, their hydrodynamic interaction can be represented with
a moderate number of points, but as they approach one another, the
hydrodynamic interaction strengthens and this requires a finer spatial
resolution.

%% file: aver.tex
In many applications, we need to compute local averages of pressure and stress.  The pressure $p^{S}$ and stress $T^{S}$ of the single-layer potential $\uu(\xx) = \SS[\ff](\xx)$ are 
\begin{align*}
  p^{S}(\xx) = \sum_{q=1}^{M} p_{q}^{S}(\xx) = 
    \frac{1}{2\pi} \sum_{q=1}^{M}\int_{\gamma_{q}}
    \frac{\rr \cdot \ff_{q}}{\rho^{2}} ds,
\end{align*}
and
\begin{align*}
  T^{S}[\ssigma](\xx) = \sum_{q=1}^{M} T_{q}^{S}[\ssigma](\xx) = 
    \frac{1}{\pi}\sum_{q=1}^{M} \int_{\gamma_{q}} 
    \frac{\rr \cdot \ssigma}{\rho^{2}} 
    \frac{\rr \otimes \rr}{\rho^{2}}\ff_{q} ds.
\end{align*}
The pressure $p^{D}$ and stress $T^{D}$ of the double-layer potential
$\uu(\xx) = \DD[\ff](\xx)$ are
\begin{align*}
  p^{D}(\xx) = \sum_{q=1}^{M} p_{q}^{D}(\xx) =
    -\frac{1}{\pi}\sum_{q=1}^{M}\int_{\gamma_{q}}
    \frac{1}{\rho^{2}}\left(1 - 2\frac{\rr \otimes \rr}{\rho^{2}}\right)
    \nn \cdot \ff_{q} ds,
\end{align*}
and 
\begin{align*}
  T^{D}[\ssigma](\xx) &= \sum_{q=1}^{M}T_{q}^{D}[\ssigma](\xx) \\
  &= \frac{1}{\pi}\sum_{q=1}^{M}\int_{\gamma_{q}}
    \left(\frac{\nn \cdot \ff_{q}}{\rho^{2}}\ssigma - 
    \frac{8}{\rho^{6}}(\rr \cdot \nn)(\rr \cdot \ff_{q})
      (\rr \cdot \ssigma)\rr \right.\\ 
    &\hspace{30pt}+\left.\frac{\rr \cdot \nn}{\rho^{4}}
      (\rr \otimes \ff_{q} + \ff_{q} \otimes \rr)\ssigma +
    \frac{\rr \cdot \ff_{q}}{\rho^{4}}
      (\rr \otimes \nn + \nn \otimes \rr)\ssigma \right)ds.
\end{align*}
With these expressions, we can evaluate the pressure and stress tensor
for bounded and unbounded flows with viscosity contrast.  Moreover,
computing the pressure of the single-layer potential can be
accelerated with the standard Laplace FMM.

To evaluate the pressure and stress tensor near a vesicle, we again
require near-singular integration.  Since we require boundary values
(see $\xx_{0}$ in the left plot of Figure~\ref{f:nearsing:conv}),
formulas for the jumps in pressure and stress as a target point
approaches a vesicle are required.  The pressures $p_{q}^{S}$ and
$p_{q}^{D}$ satisfy
\begin{align*}
  \lim_{\substack{\xx \rightarrow \xx_{0} \\ \xx \in \omega_{q}}} p^{S}_{q}(\xx) &= 
   \frac{\ff_{0} \cdot \nn_{0}}{2} + p^{S}_{q}(\xx_{0}), 
  &&\lim_{\substack{\xx \rightarrow \xx_{0} \\ \xx \notin \omega_{q}}} p^{S}_{q}(\xx) = 
   -\frac{\ff_{0} \cdot \nn_{0}}{2} + p^{S}_{q}(\xx_{0}), \\
  \lim_{\substack{\xx \rightarrow \xx_{0} \\ \xx \in \omega_{q}}} p^{D}_{q}(\xx) &= 
    -\frac{\p \ff_{0}}{\p \ttau} \cdot \ttau + p^{D}_{q}(\xx_{0}), 
  &&\lim_{\substack{\xx \rightarrow \xx_{0} \\ \xx \notin \omega_{q}}} p^{D}_{q}(\xx) = 
    \frac{\p \ff_{0}}{\p \ttau} \cdot \ttau + p^{D}_{q}(\xx_{0}),
\end{align*}
where $\xx_{0} \in \gamma$ and $\ff_{0} = \ff(\xx_{0})$.  The jumps for
$p_{q}^{S}$ are proved in Appendix~\ref{A:AppendixB} and the jumps for
$p_{q}^{D}$ are proved in~\cite{ying-biros-zorin06}.  To compute
$p_{q}^{S}$ and $p_{q}^{D}$ at $\xx_{0}$, we use odd-even integration
which has spectral accuracy if the singularity of the integrand at
$\xx=\xx_{0}$ is no stronger than $(\xx -
\xx_{0})^{-1}$~\cite{sid:isr}.  Odd-even integration can be directly
applied to evaluate $p_{q}^{S}(\xx_{0})$.  However, as $\xx \rightarrow
\xx_{0}$, the singularity of $p_{q}^{D}(\xx_{0})$ is of the order
$(\xx - \xx_{0})^{-2}$; therefore, we have to interpret the integral in
the principal value sense.  Since a constant hydrodynamic density $\ff$
results in the pressure vanishing~\cite{kress,pozrikidis1992}, we can
reduce the strength of the singularity to $(\xx - \xx_{0})^{-1}$ by
first subtracting $\ff(\xx_{0})$ and then evaluating
$p^{D}_{q}(\xx_{0})$.

We do a convergence study on $p^{S}$ and $p^{D}$ exterior to the unit
circle with the hydrodynamic density $\ff = [\cos(\theta)
\cos(\theta)]^{T}$.  The exact pressures are calculated using the
Residue Theorem.  We check the relative maximum errors along the
vertical line $x=1.01$ (Table~\ref{t:pressure:conv}).  Since
$p^{S},p^{D} \in C^{\infty}$, we achieve super-algebraic convergence
when near-singular integration is not required.  However, near-singular
integration introduces an error described in
Appendix~\ref{A:AppendixA}.  Table~\ref{t:pressure:conv} indicates a
$5^{th}$-order convergence rate which is consistent with the
$6^{th}$-degree Lagrange interpolant used to interpolate values on
$\gamma$ to $\xx_{0}$.  Figure~\ref{f:pressure:figure} shows contour
plots of the pressure for an extensional flow with no viscosity
contrast, as well as the pressure along the vertical line that passes
exactly between the two vesicles.  Not surprisingly, the pressure
between the two vesicles is increasing as the inter-vesicle distance
decreases.
\begin{table}[htps]
\begin{centering}
\begin{tabular}{l|cccccc}
& $N=32$ & $N=64$ & $N=128$ & $N=256$ & $N=512$ & $N=1024$ \\
\hline
Single-layer Potential & 
  $8.22e-04$ & $4.05e-05$ & $1.23e-06$ & $2.53e-08$ & $2.03e-10$ & $1.83e-14$ \\ 
Double-layer Potential & 
  $8.22e-04$ & $4.05e-05$ & $1.23e-06$ & $2.53e-08$ & $2.03e-10$ & $2.78e-13$ 
\end{tabular}
\mcaption{The maximum relative errors in the calculation of the
pressures using near-singular integration along the line $x=1.01$.  The
exact pressures are computed analytically using the Residue Theorem,
and they differ by a multiplicative constant which explains the
proximity of the errors.  The slope of the line of best fit on a
log-log scale is 5.5.  For $N=1024$, the near zone $\Omega_{0}$ is
empty.  This explains the sharp drop in the errors.}{t:pressure:conv}
\end{centering}
\end{table}

\begin{figure}[htps]
\centering
$
\begin{array}{ccccc}
  \includegraphics[trim=1.2cm 7cm 2cm 6cm,clip=true,scale = 0.15]{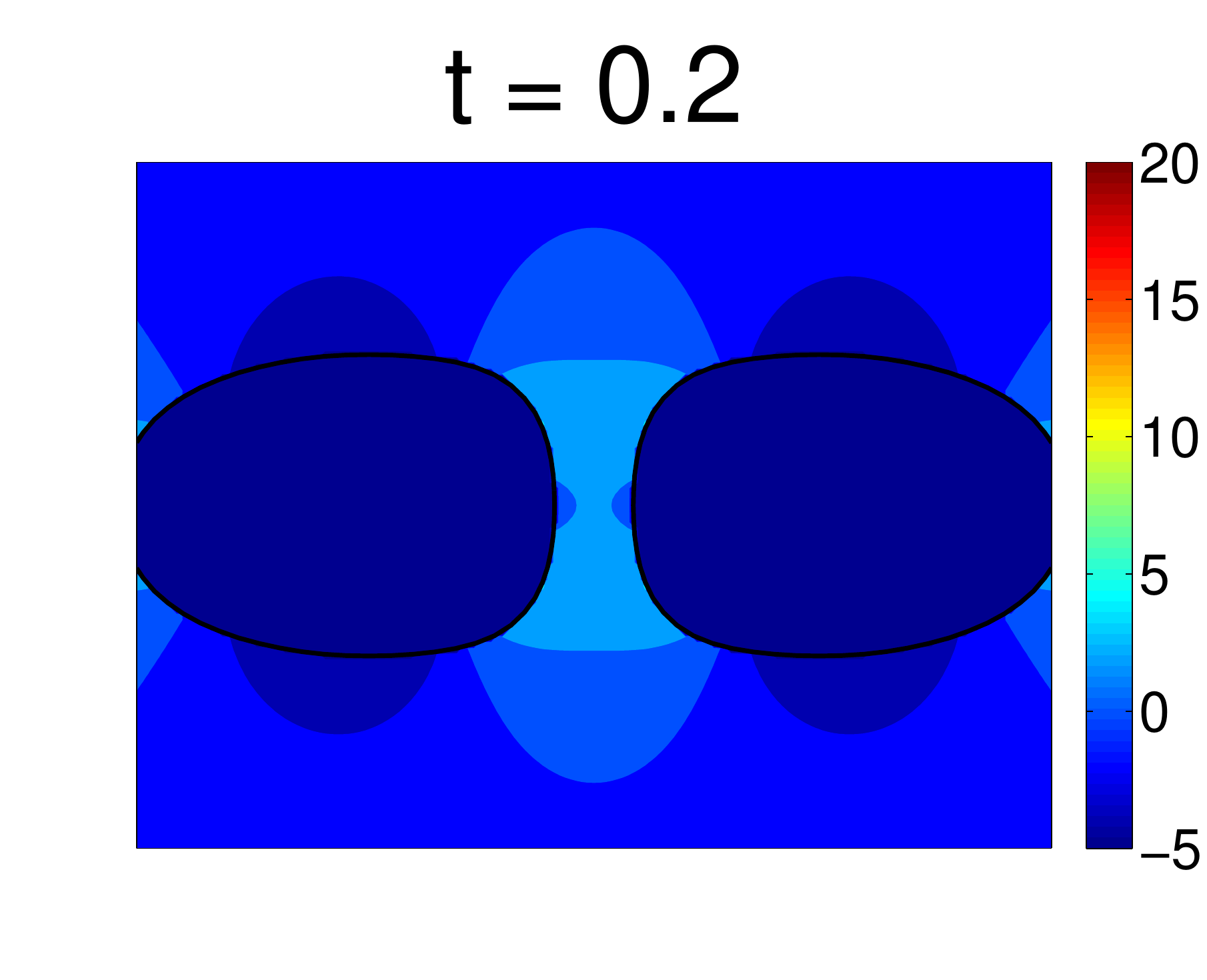} &
  \includegraphics[trim=1.2cm 7cm 2cm 6cm,clip=true,scale = 0.15]{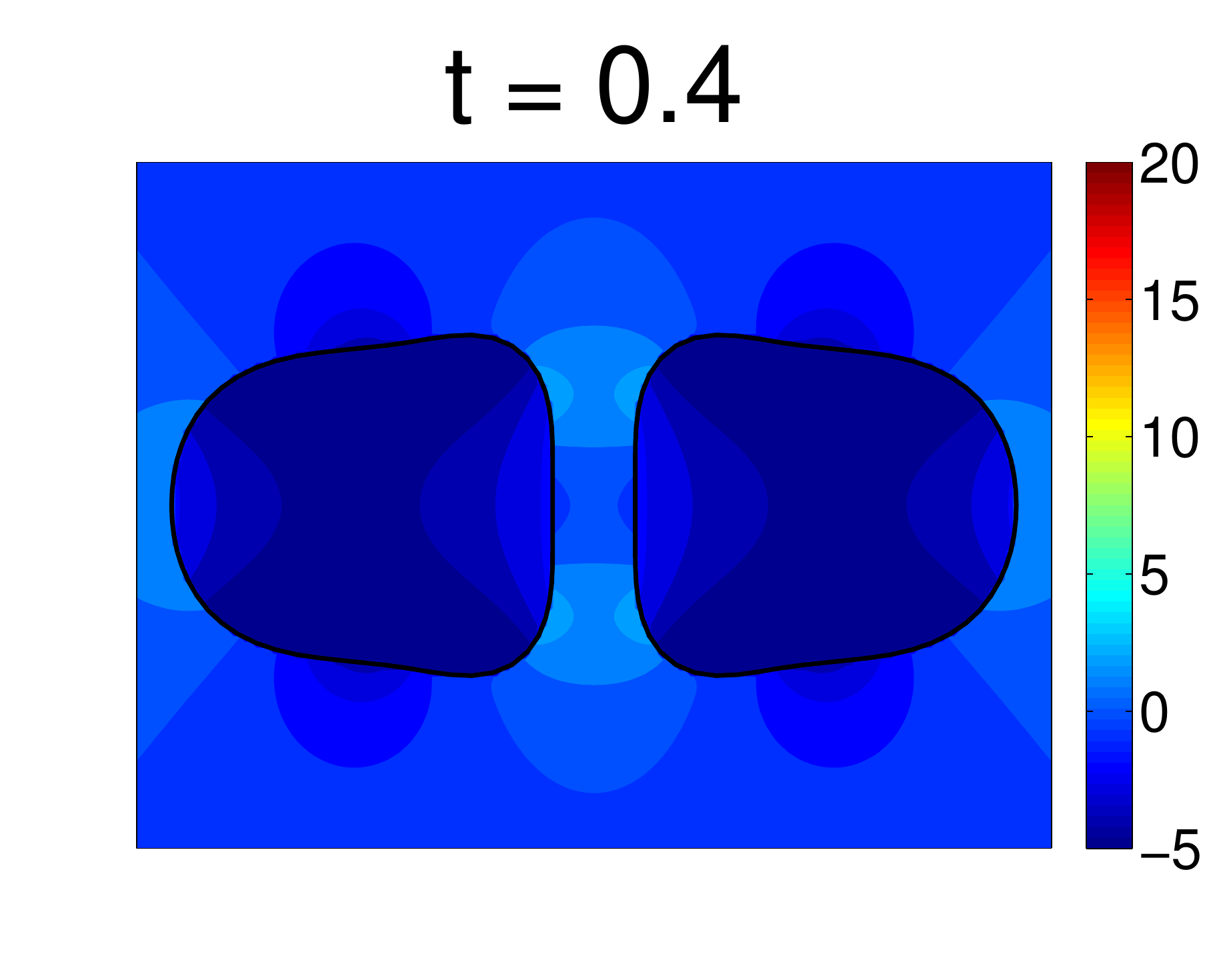} &
  \includegraphics[trim=1.2cm 7cm 2cm 6cm,clip=true,scale = 0.15]{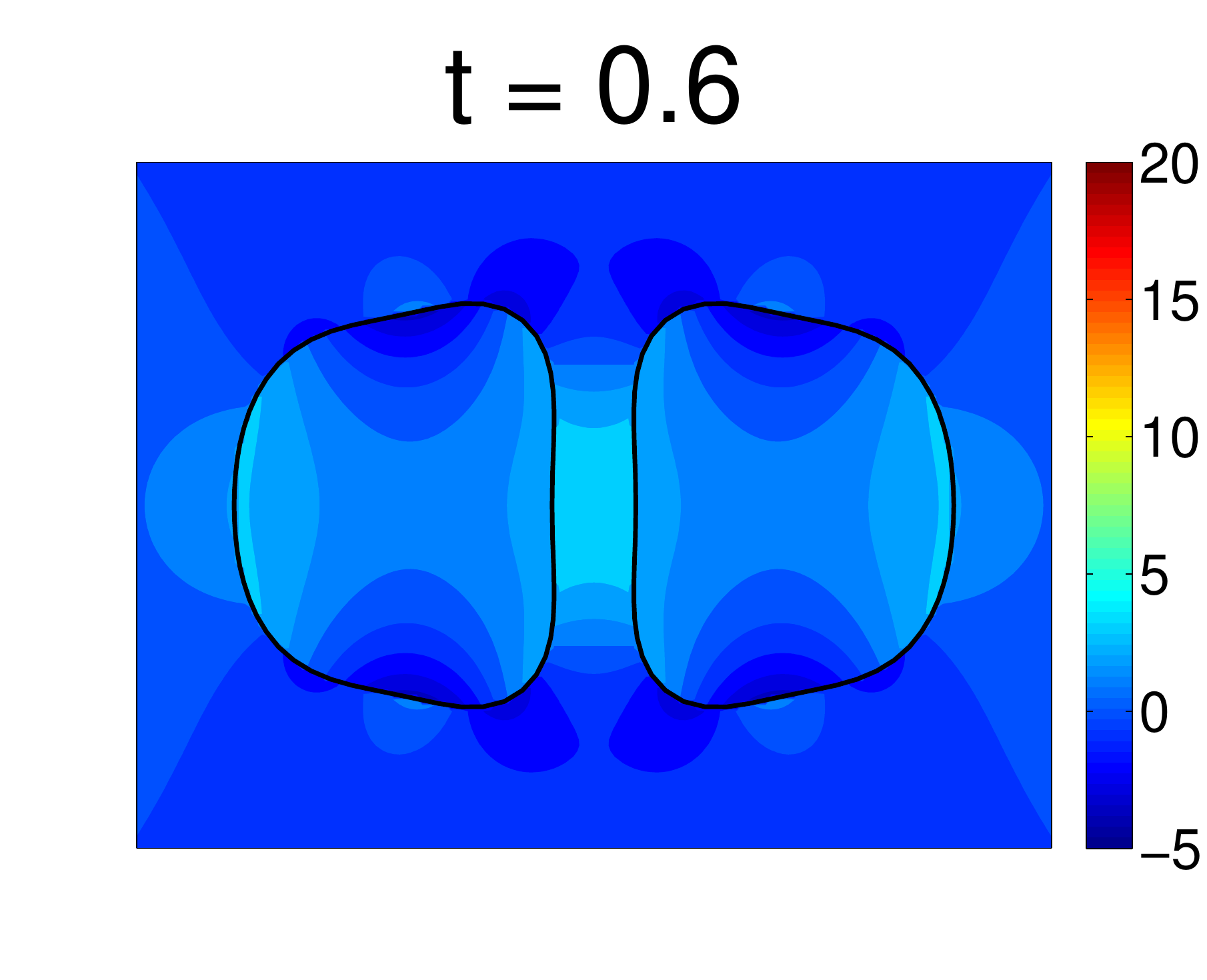} &
  \includegraphics[trim=1.2cm 7cm 2cm 6cm,clip=true,scale = 0.15]{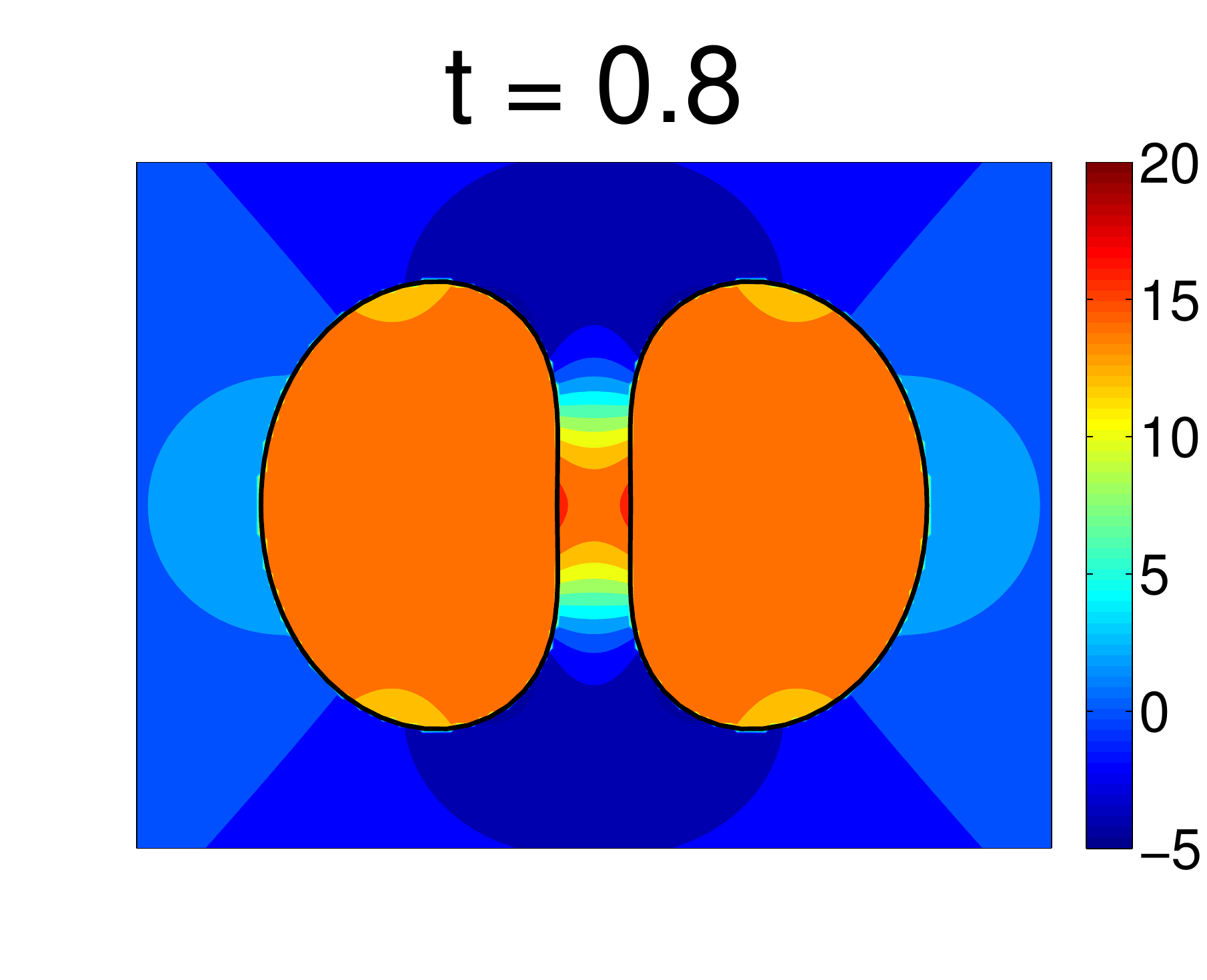} &
  \includegraphics[trim=1.2cm 7cm 2cm 6cm,clip=true,scale = 0.15]{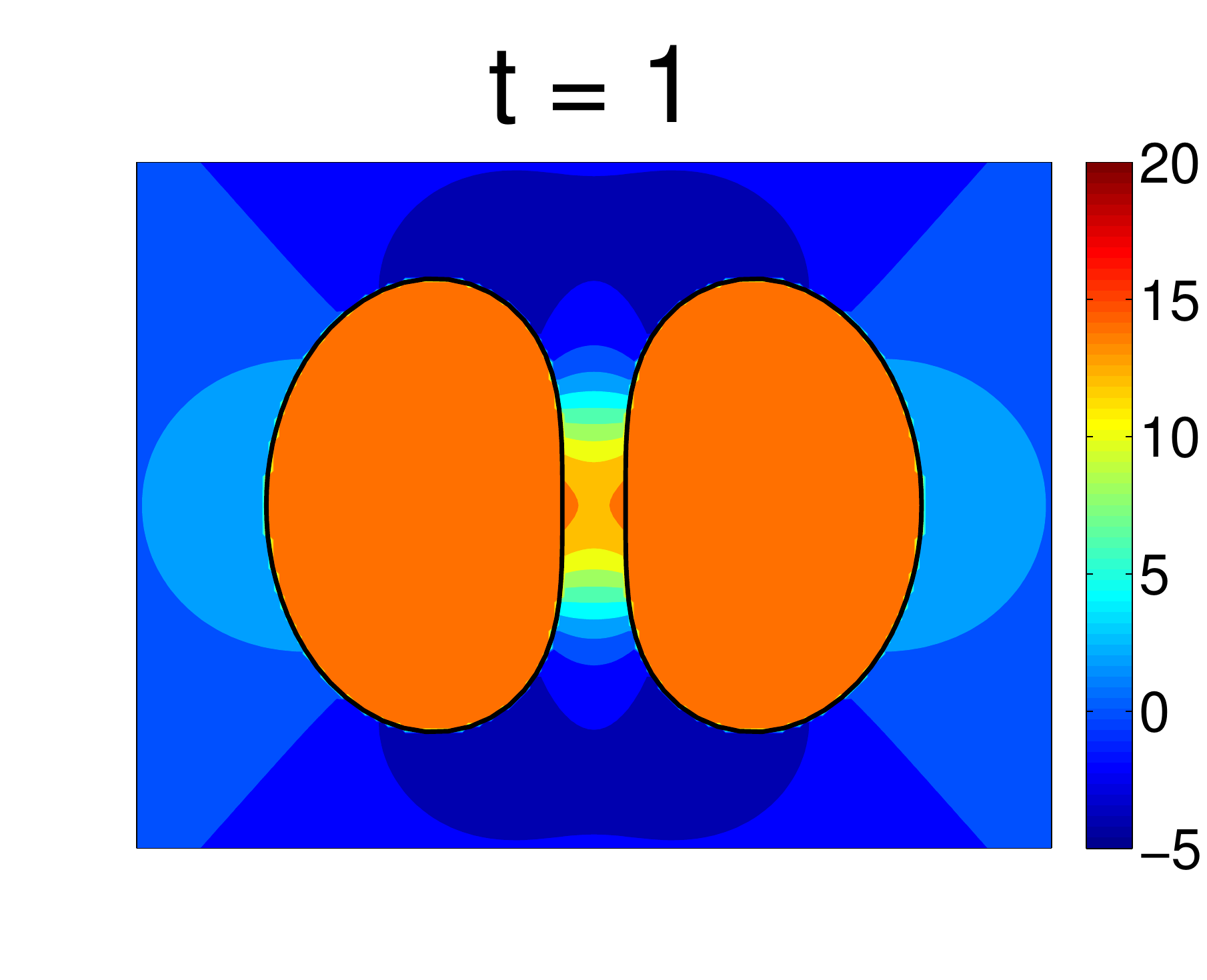} \\ 
  \includegraphics[trim=1.2cm 7cm 2cm 6cm,clip=true,scale = 0.15]{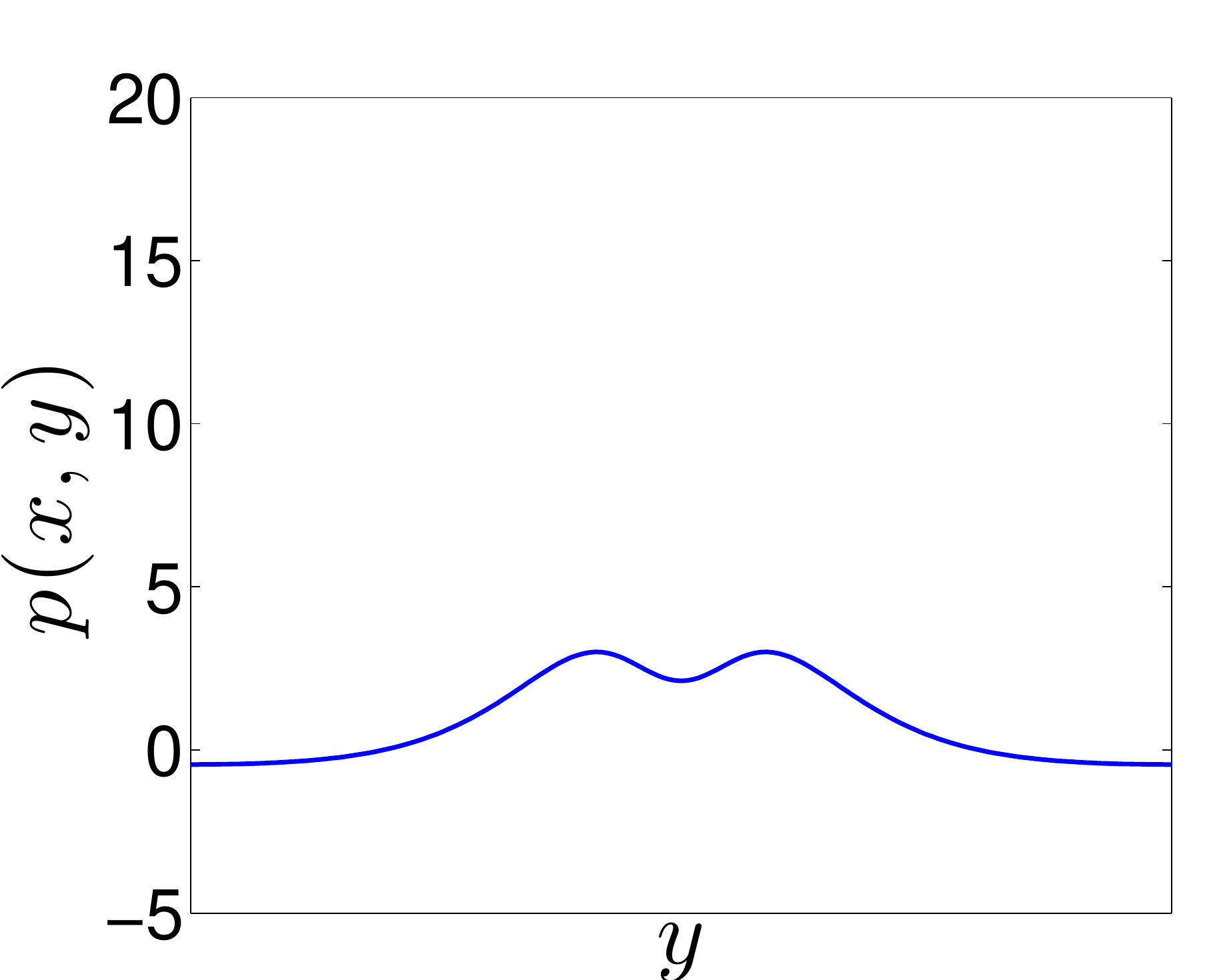} &
  \includegraphics[trim=1.2cm 7cm 2cm 6cm,clip=true,scale = 0.15]{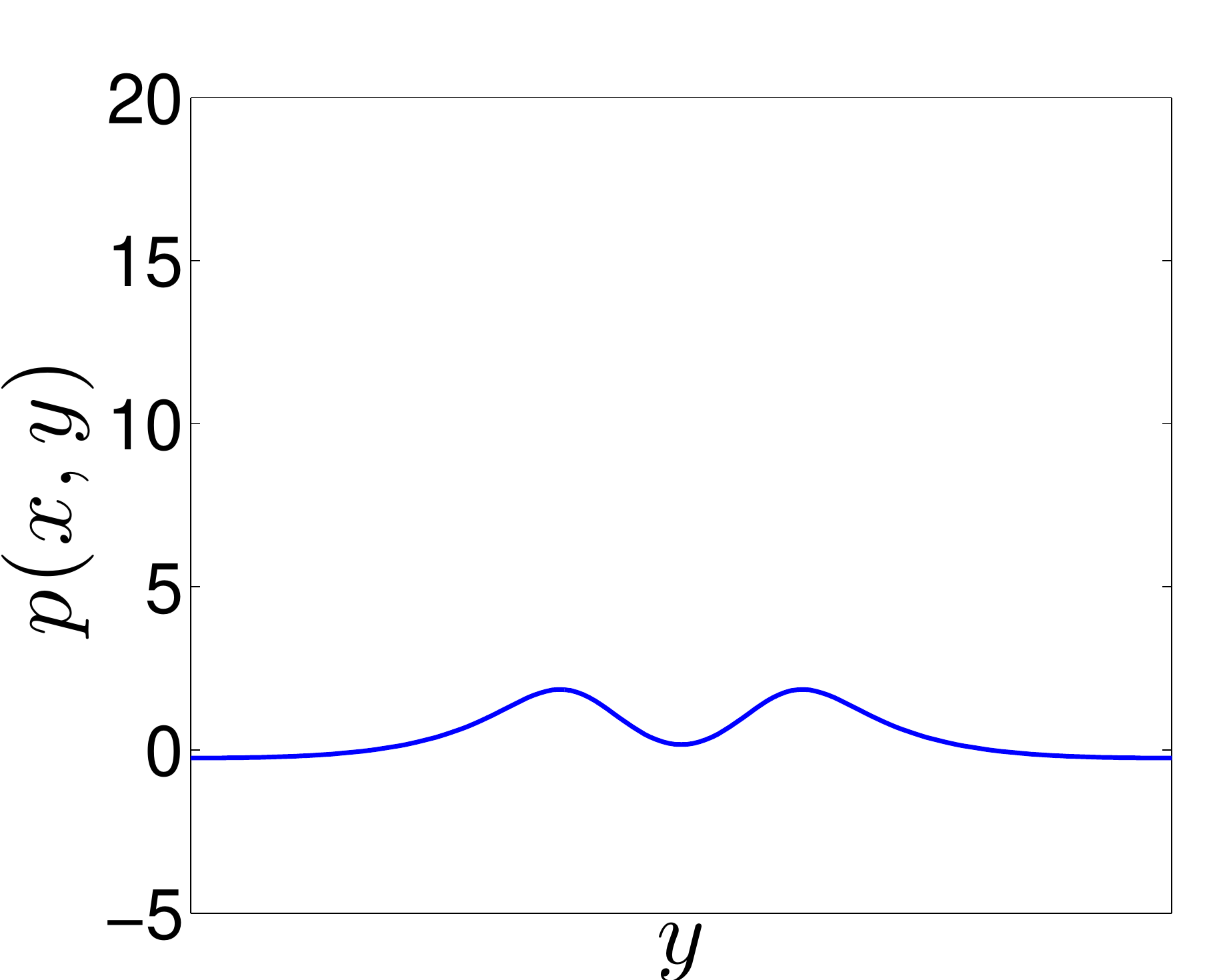} &
  \includegraphics[trim=1.2cm 7cm 2cm 6cm,clip=true,scale = 0.15]{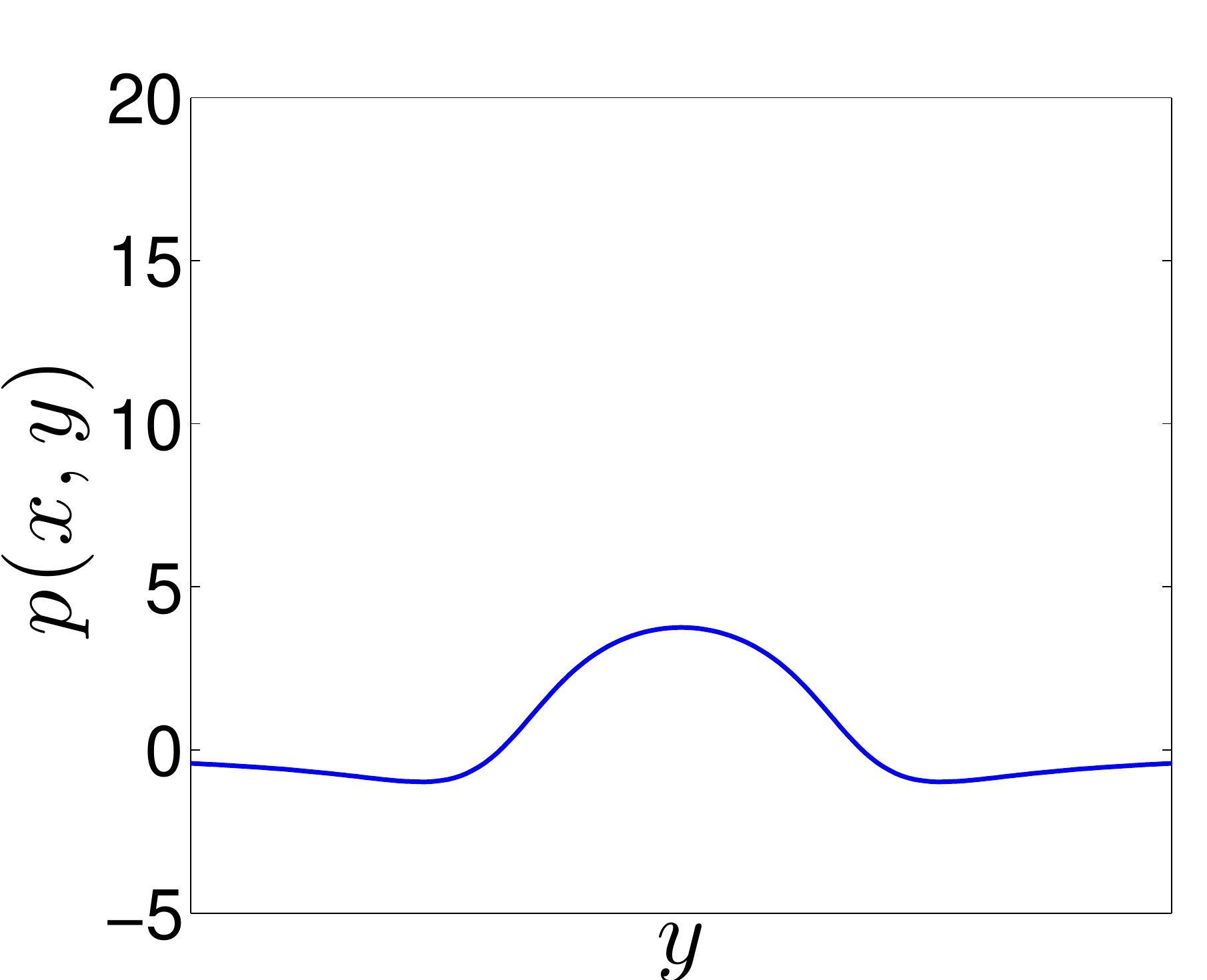} &
  \includegraphics[trim=1.2cm 7cm 2cm 6cm,clip=true,scale = 0.15]{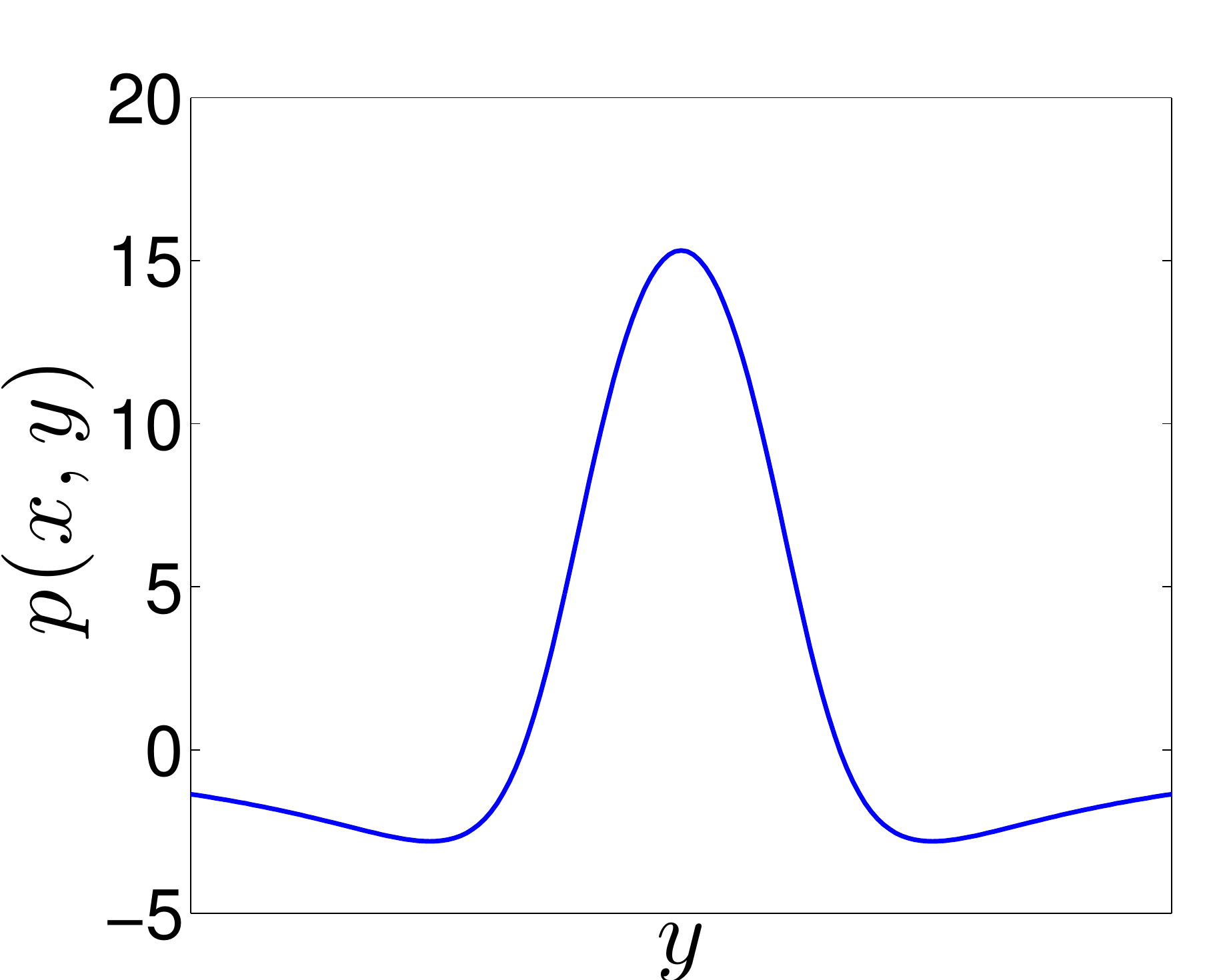} &
  \includegraphics[trim=1.2cm 7cm 2cm 6cm,clip=true,scale = 0.15]{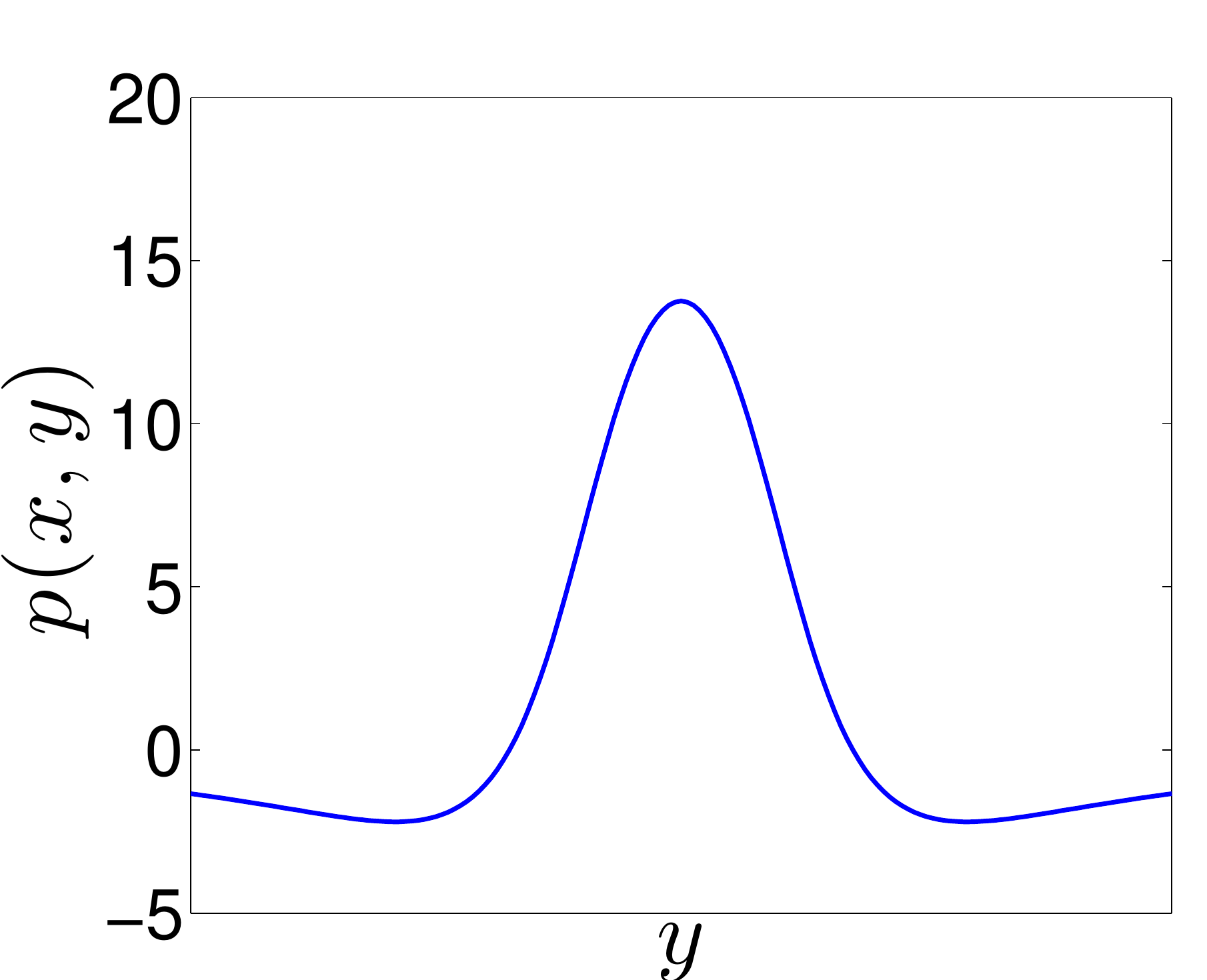} \\ \\ 
  \includegraphics[trim=1.2cm 7cm 2cm 6cm,clip=true,scale = 0.15]{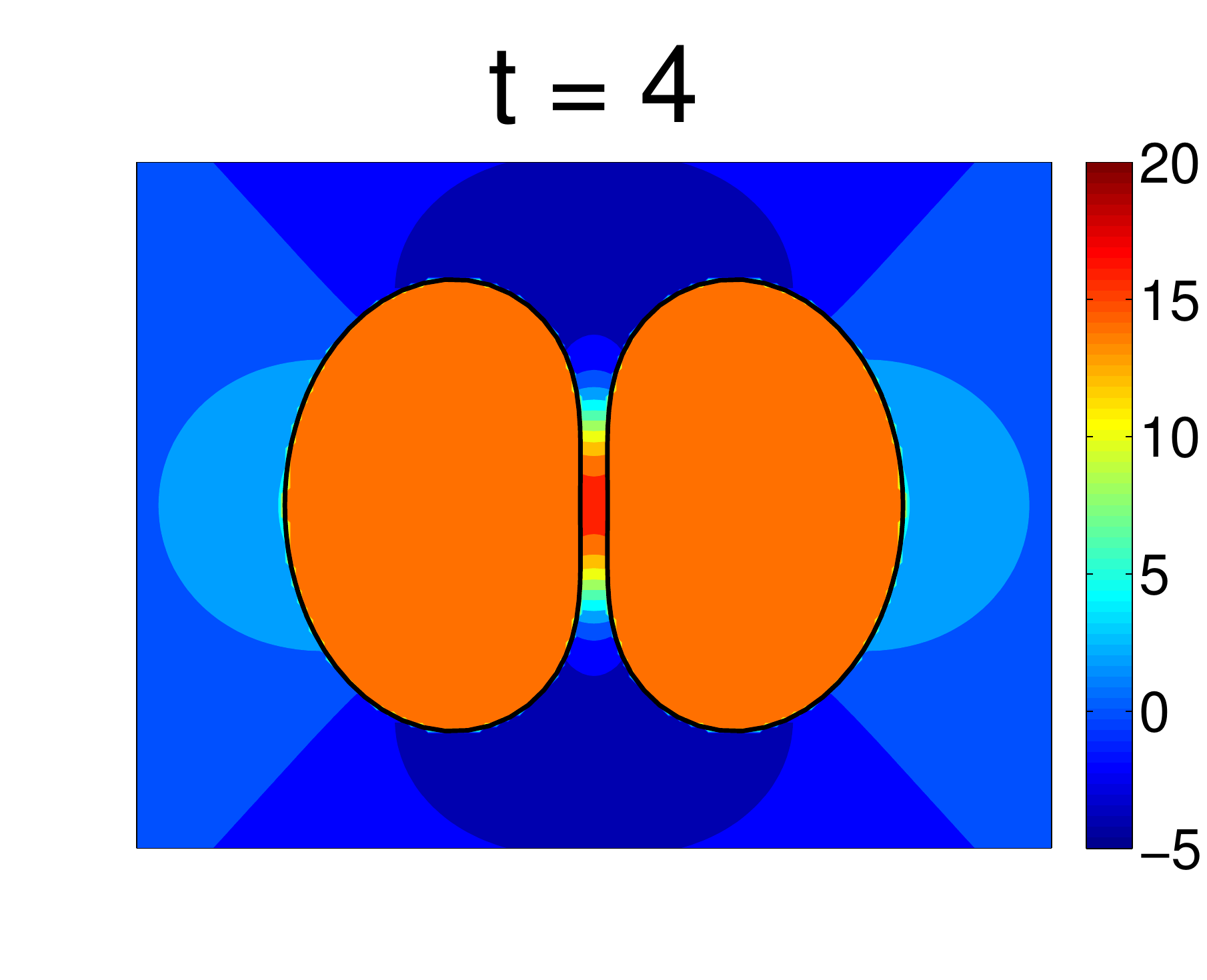} &
  \includegraphics[trim=1.2cm 7cm 2cm 6cm,clip=true,scale = 0.15]{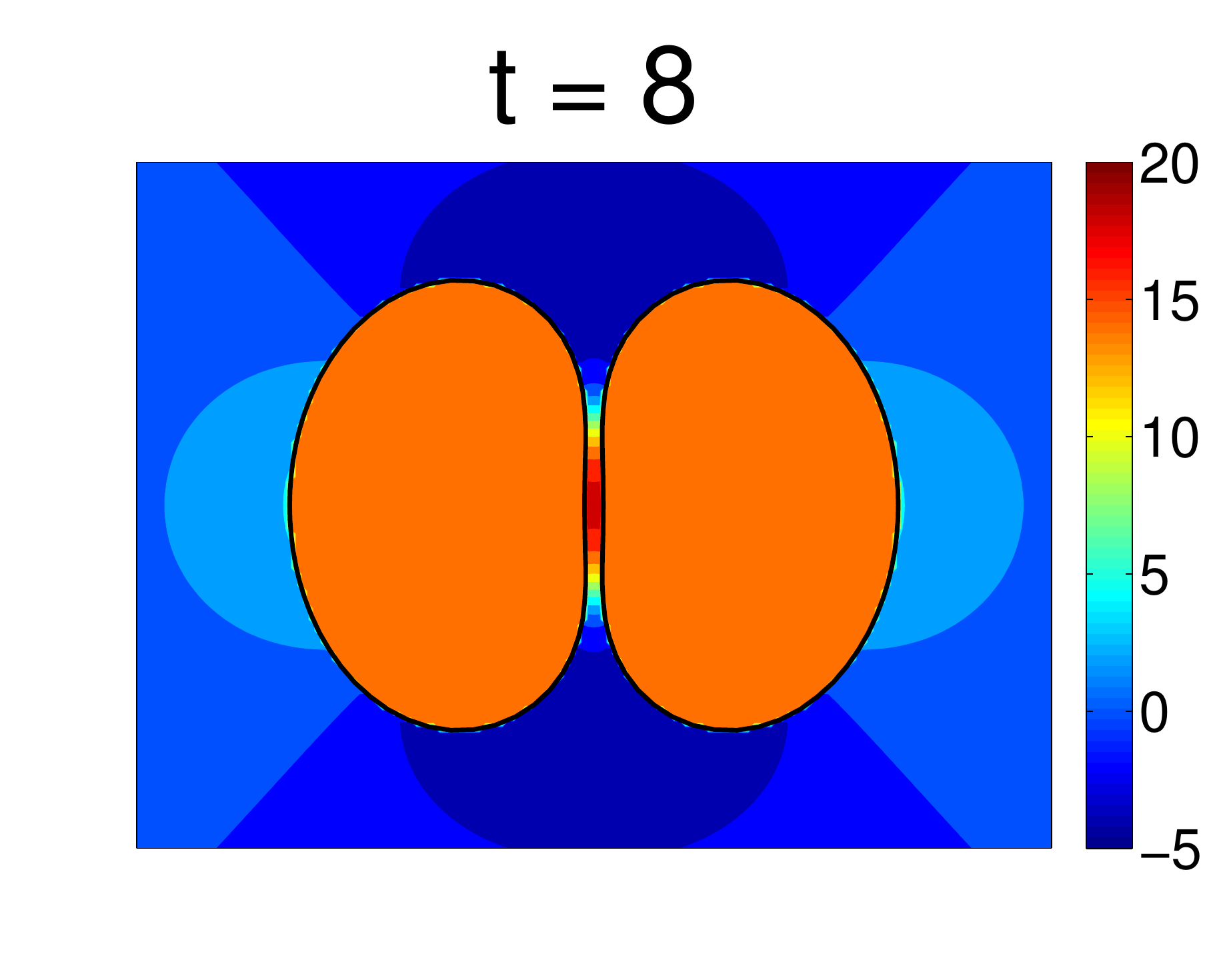} &
  \includegraphics[trim=1.2cm 7cm 2cm 6cm,clip=true,scale = 0.15]{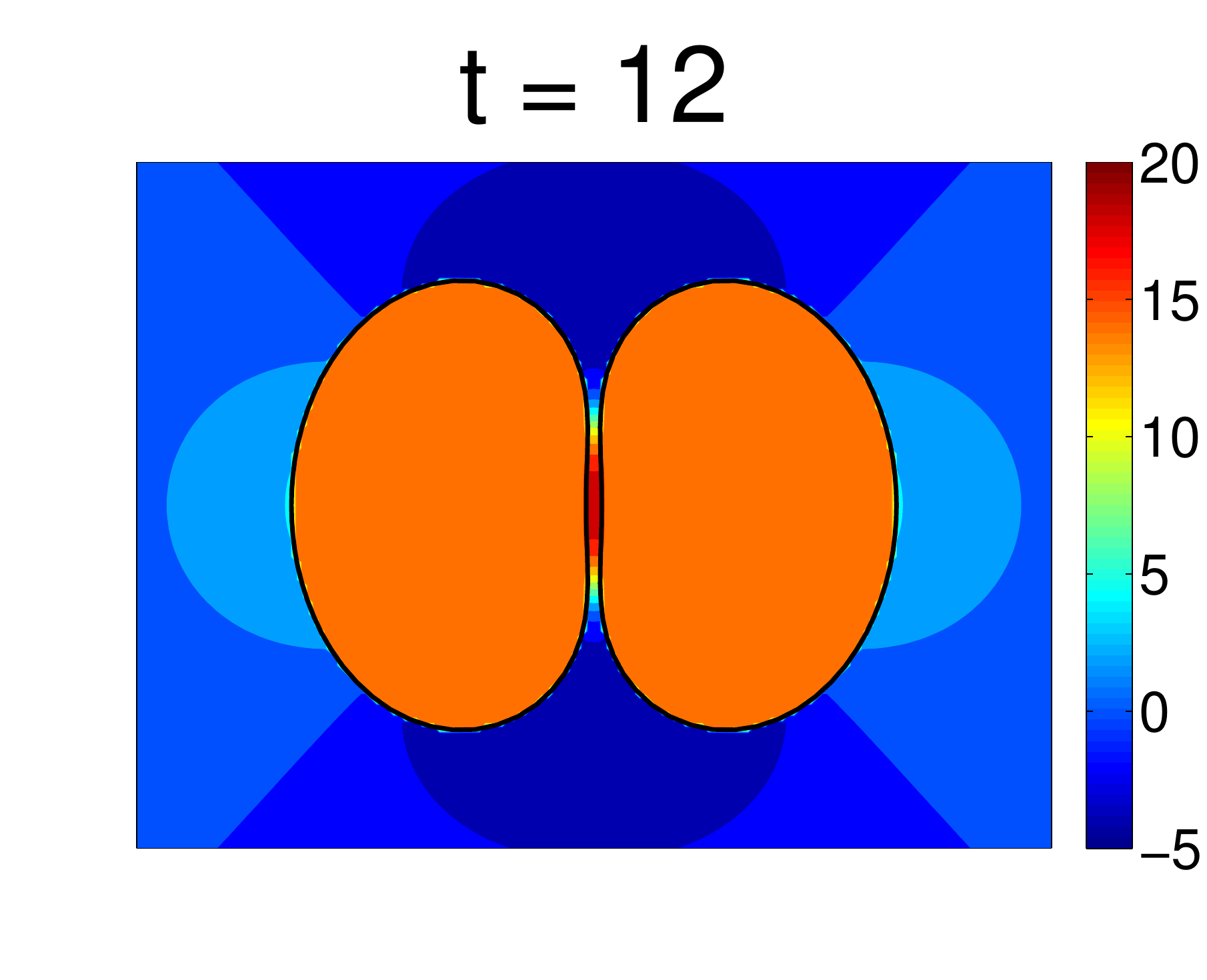} &
  \includegraphics[trim=1.2cm 7cm 2cm 6cm,clip=true,scale = 0.15]{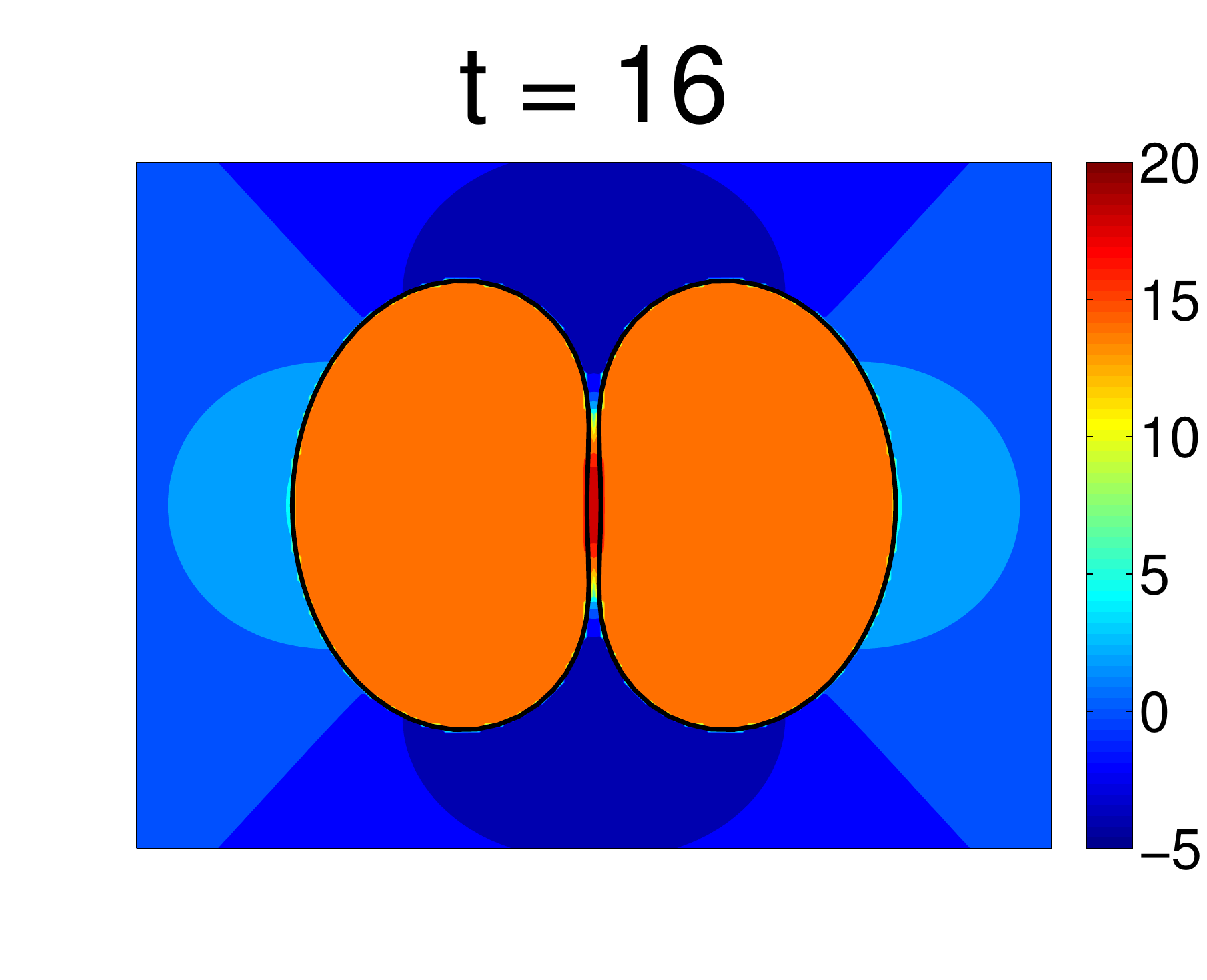} &
  \includegraphics[trim=1.2cm 7cm 2cm 6cm,clip=true,scale = 0.15]{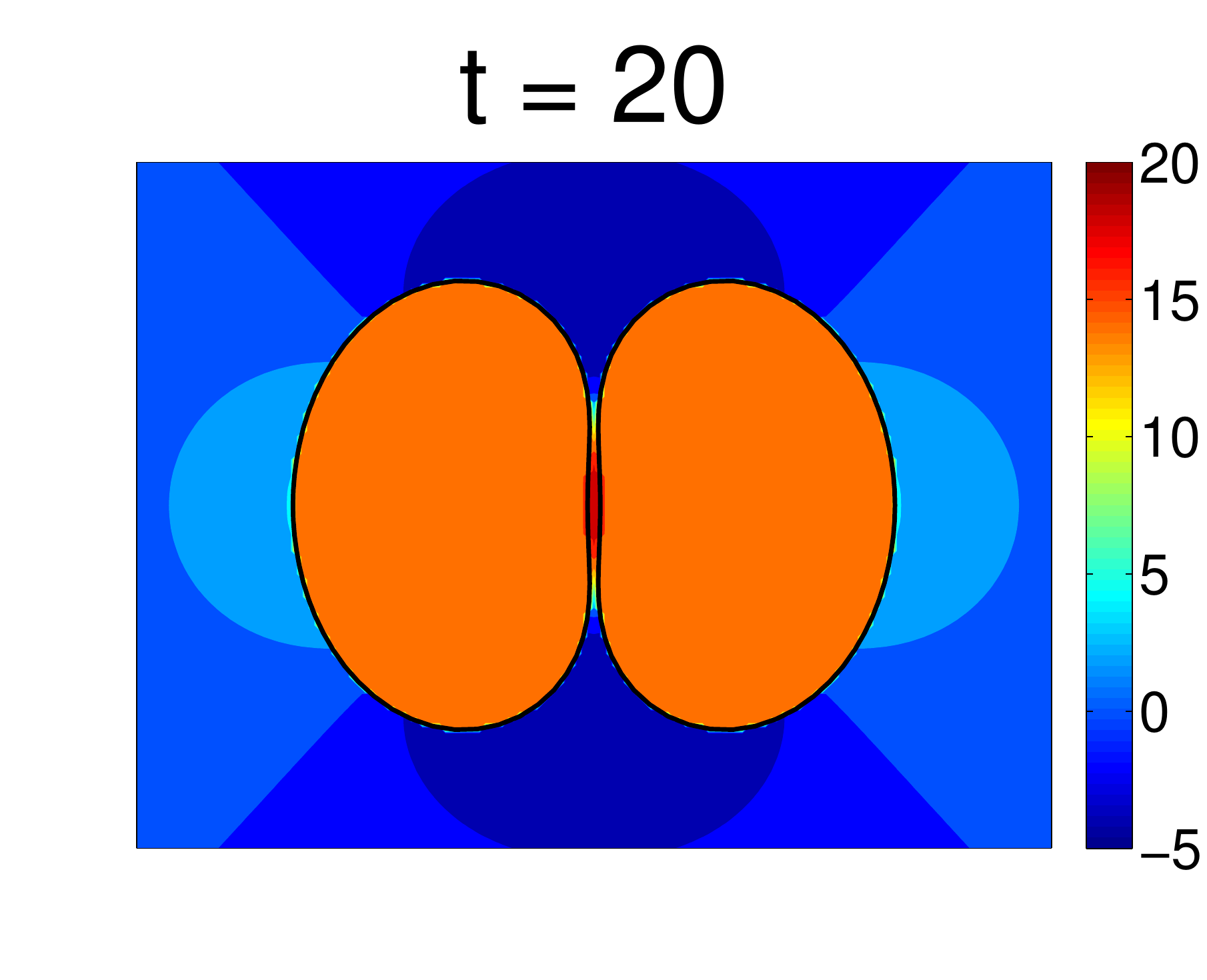} \\ 
  \includegraphics[trim=1.2cm 7cm 2cm 6cm,clip=true,scale = 0.15]{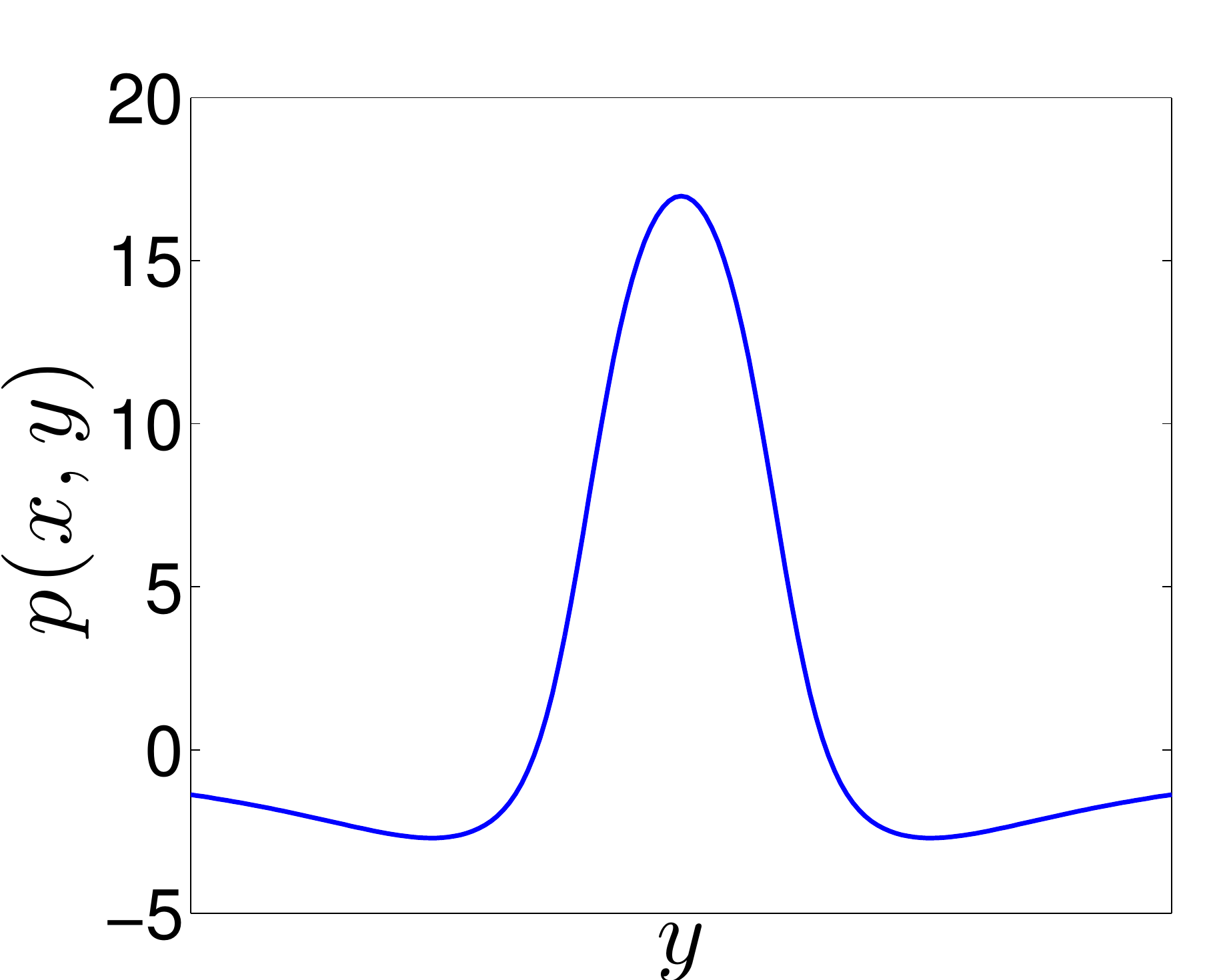} &
  \includegraphics[trim=1.2cm 7cm 2cm 6cm,clip=true,scale = 0.15]{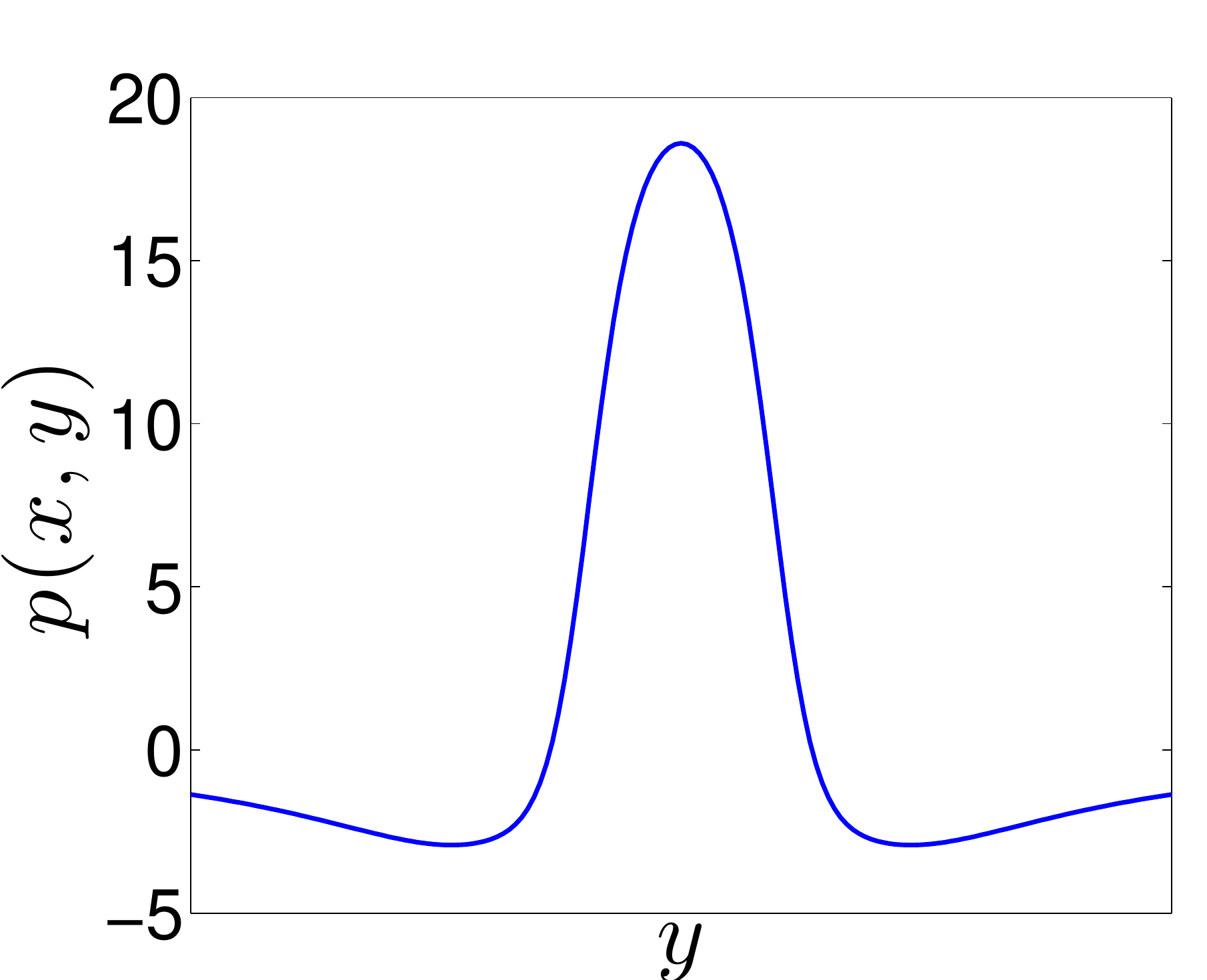} &
  \includegraphics[trim=1.2cm 7cm 2cm 6cm,clip=true,scale = 0.15]{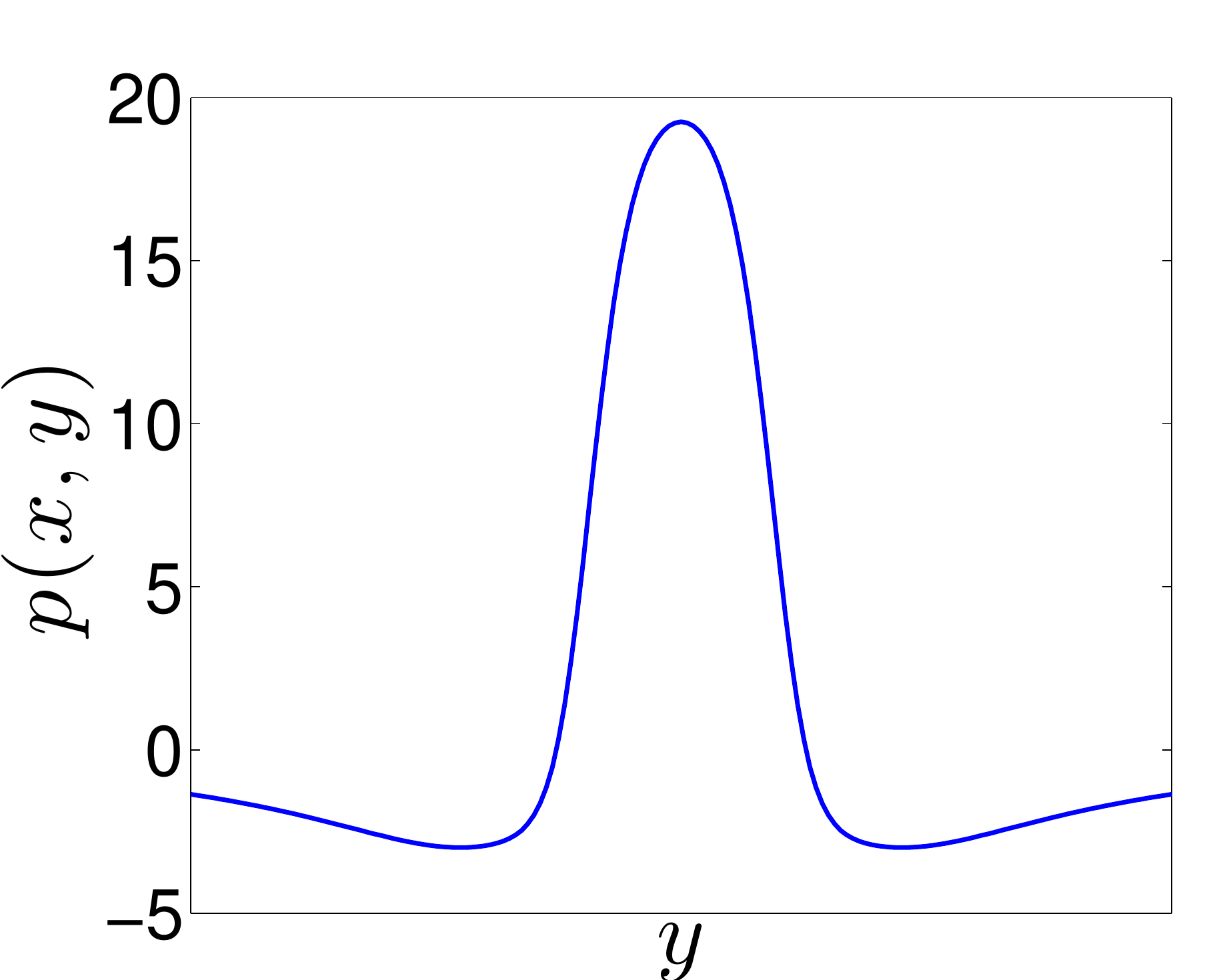} &
  \includegraphics[trim=1.2cm 7cm 2cm 6cm,clip=true,scale = 0.15]{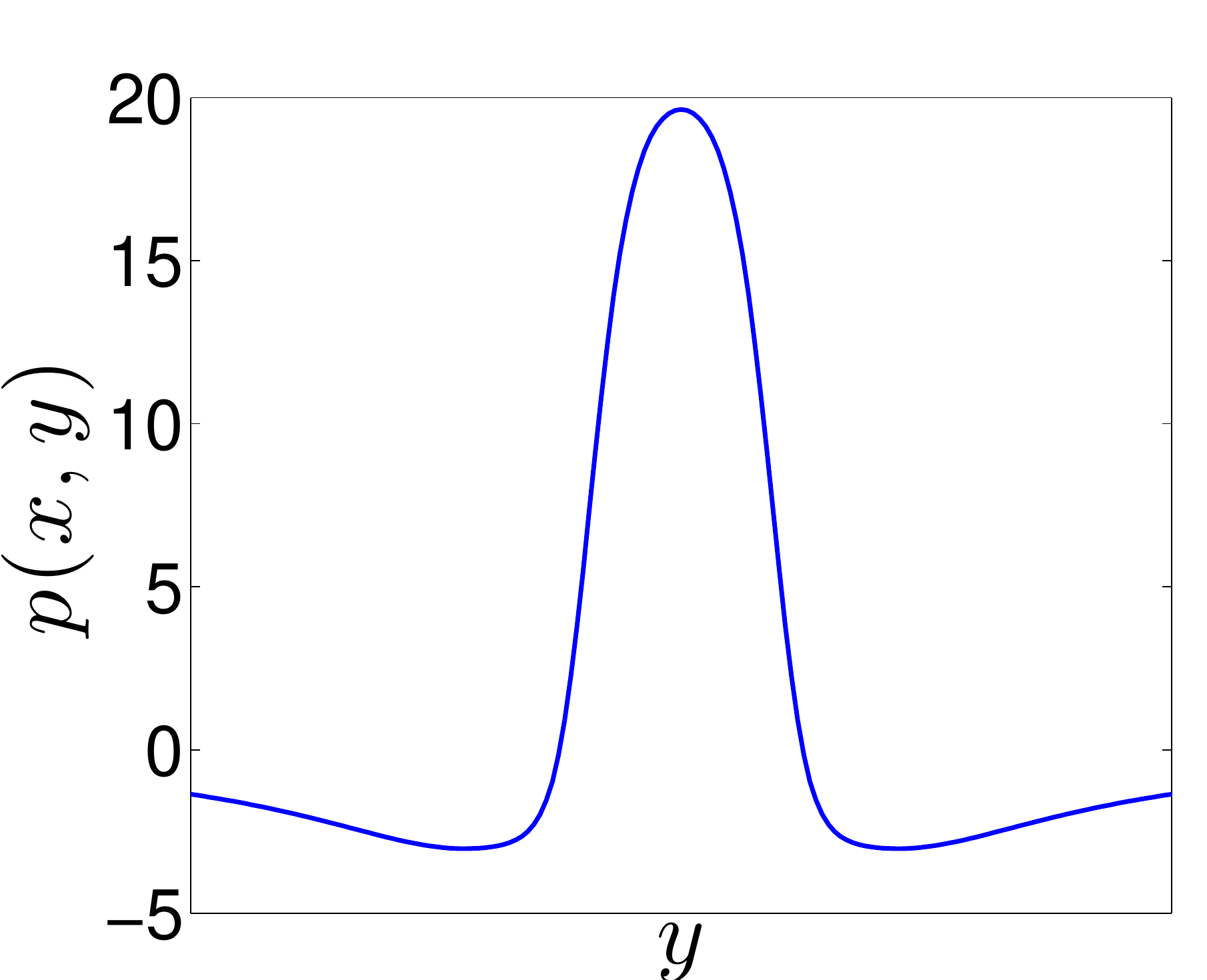} &
  \includegraphics[trim=1.2cm 7cm 2cm 6cm,clip=true,scale = 0.15]{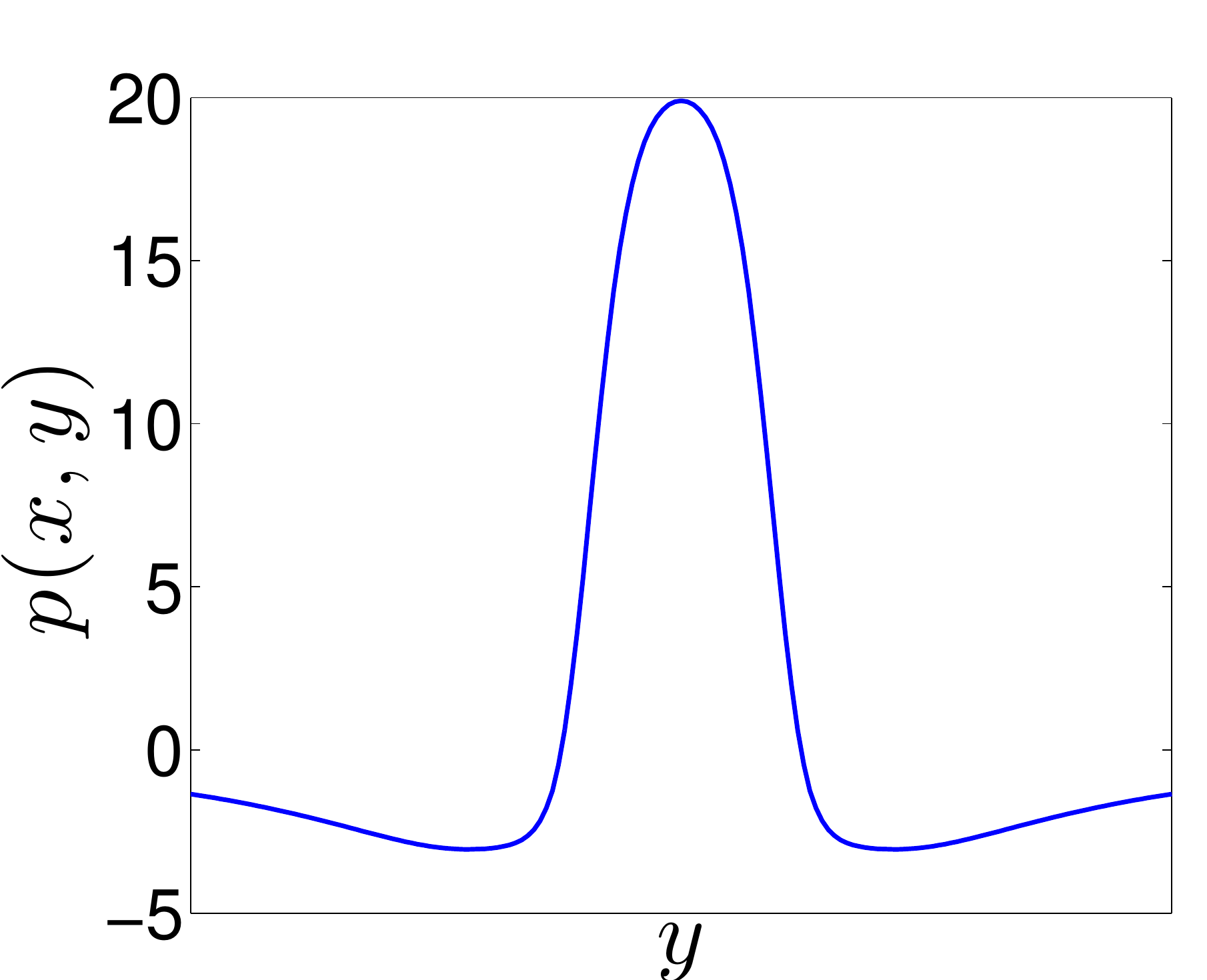} \\ 
\end{array}
$
\mcaption{Contour plots show the pressure in and around the vesicles in
an extensional flow with no viscosity contrast.  The lower plots are
the pressure along the vertical line that passes vertically through the
midpoint of the two vesicles.  We see that the maximum pressure
increases as a function of time and that the pressure inside each
vesicle approaches a constant value.}{f:pressure:figure}
\end{figure}

We now compute the limiting values of the stress tensors
$T^{S}[\ssigma]$ and $T^{D}[\ssigma]$.  The limiting values of
$T^{S}_{q}[\ssigma]$ are
\begin{align}
  \lim_{\substack{\xx \rightarrow \xx_{0} \\ \xx \in \omega_{q}}} 
    T^{S}_{q}[\ssigma](\xx) 
  &=-\frac{1}{2}(\nn_{0} \otimes \ff_{0})\ssigma_{0} +
    \frac{1}{2}\left( \ttau \otimes \left[
    \begin{array}{cc}
      2\tau_{x}\tau_{y} & \tau_{y}^{2} - \tau_{x}^{2} \\
      \tau_{y}^{2} - \tau_{x}^{2} & -2\tau_{x}\tau_{y}
    \end{array}
    \right]\ff_{0} \right) \ssigma_{0} + T^{S}_{q}[\ssigma](\xx_{0}), 
    \label{e:stress:jump:int} \\
  \lim_{\substack{\xx \rightarrow \xx_{0} \\ \xx \notin \omega_{q}}} 
    T^{S}_{q}[\ssigma](\xx) 
  &=\frac{1}{2}(\nn_{0} \otimes \ff_{0})\ssigma_{0} -
    \frac{1}{2}\left( \ttau \otimes \left[
    \begin{array}{cc}
      2\tau_{x}\tau_{y} & \tau_{y}^{2} - \tau_{x}^{2} \\
      \tau_{y}^{2} - \tau_{x}^{2} & -2\tau_{x}\tau_{y}
    \end{array}
    \right]\ff_{0} \right) \ssigma_{0} + T^{S}_{q}[\ssigma](\xx_{0}).
    \label{e:stress:jump:ext}
\end{align}
The jumps are proved in Appendix~\ref{A:AppendixB}.  Since the
singularity at $\xx-\xx_{0}$ is of the order $(\xx - \xx_{0})^{-1}$, we
use odd-even integration to compute the stress on $\gamma_{q}$.  The
jumps in the tensor of the double-layer potential are
\begin{align*}
  \lim_{\substack{\xx \rightarrow \xx_{0} \\ \xx \in \omega_{q}}} 
    T^{D}_{q}[\ssigma](\xx) &= -\frac{\p \ff_{0}}{\p \ttau} \cdot \ttau
    \left(I + \left[
    \begin{array}{cc}
      \tau_{x}^{2}-\tau_{y}^{2} & 2\tau_{x}\tau_{y} \\
      2\tau_{x}\tau_{y} & -\tau_{x}^{2} + \tau_{y}^{2}
    \end{array}
    \right] \right) \ssigma_{0} + T^{D}_{q}[\ssigma](\xx_{0}),
  \\
  \lim_{\substack{\xx \rightarrow \xx_{0} \\ \xx \notin \omega_{q}}} 
    T^{D}_{q}[\ssigma](\xx) &= \frac{\p \ff_{0}}{\p \ttau} \cdot \ttau
    \left(I + \left[
    \begin{array}{cc}
      \tau_{x}^{2}-\tau_{y}^{2} & 2\tau_{x}\tau_{y} \\
      2\tau_{x}\tau_{y} & -\tau_{x}^{2} + \tau_{y}^{2}
    \end{array}
    \right] \right) \ssigma_{0} + T^{D}_{q}[\ssigma](\xx_{0}).
\end{align*}
These jumps are calculated using the jumps $\llbracket p^{D}
\rrbracket$, $\llbracket \uu^{D} \rrbracket$, $\llbracket T^{D}\nn
\rrbracket$, and $\llbracket \grad \cdot \uu^{D} \rrbracket$.  We have
already computed $\llbracket p^{D} \rrbracket$, and the other three
jumps can be found in Appendix B of~\cite{ying-biros-zorin06}.
However, we have to use the same technique we used to compute
$p^{D}(\xx_{0})$ to reduce the order of the singularity of the integral
from $(\xx-\xx_{0})^{-2}$ to $(\xx-\xx_{0})^{-1}$.  Since a constant
hydrodynamic density $\ff$ corresponds to a vanishing stress tensor
$T^{D}$, this modification guarantees that odd-even integration will
converge to the correct value.

We check the convergence rate for the three components of the stress
tensors (the $(1,2)$ and $(2,1)$ entries are identical).  We use the
same vesicle and hydrodynamic density as in
Table~\ref{t:pressure:conv}.  The maximum relative errors are in
Table~\ref{t:stress:conv}.  Since $T^{S},T^{D} \in C^{\infty}$, the
errors in Table~\ref{t:stress:conv} are determined by the order of the
near-singular integration algorithm.  We obtain the expected
$5^{th}$-order convergence.
\begin{table}[htps]
\begin{centering}
\begin{tabular}{lc|cccccc}
& & $N=32$ & $N=64$ & $N=128$ & $N=256$ & $N=512$ & $N=1024$ \\
\hline
Single- & $(1,1)$ & $8.33e-03$ & $6.32e-04$ & $2.40e-05$ & $5.21e-07$ & $4.30e-09$ & $1.23e-14$ \\ 
layer & $(1,2)$ & $5.79e-03$ & $4.24e-04$ & $1.59e-05$ & $3.49e-07$ & $2.89e-09$ & $6.82e-15$ \\ 
potential & $(2,2)$ & $4.39e-03$ & $3.13e-04$ & $1.15e-05$ & $2.49e-07$ & $2.05e-09$ & $1.27e-14$ \\ 
\hline
Double- & $(1,1)$ & $4.22e-03$ & $3.34e-04$ & $1.29e-05$ & $2.81e-07$ & $2.33e-09$ & $1.12e-12$ \\ 
layer & $(1,2)$ & $1.19e-02$ & $8.98e-04$ & $3.42e-05$ & $7.54e-07$ & $6.25e-09$ & $1.75e-12$ \\ 
potential & $(2,2)$ & $9.28e-03$ & $6.80e-04$ & $2.54e-05$ & $5.50e-07$ & $4.54e-09$ & $5.36e-13$ \\ 
\end{tabular}
\mcaption{The maximum relative errors in the calculation of the stress
tensor using near-singular integration along the line $x=1.01$.  The
exact stress is computed analytically using the Residue Theorem.  The
convergence rates are all about $5.2$.  For $N=1024$, the near zone
$\Omega_{0}$ is empty.  This explains the sharp drop in the
errors.}{t:stress:conv}
\end{centering}
\end{table}

%% file: results.tex
We discuss the behaviour of the different time integrators for several
bounded and unbounded vesicle suspensions.  Because of the
incompressibility of the fluid and the inextensibility of the vesicles,
the area enclosed by a vesicle and its total length remain fixed.
However, due to numerical errors, these quantities change throughout
simulations.  Penalty methods are often used to enforce maintain the
vesicles length, but these techniques change the physics of the problem.
In our formulation, no such penalty method is used and the
inextensibility constraint is enforced algebraically
using~\eqref{e:inextensibility}.  Even though the length and area of the
vesicles are low-order moments of the shape, we have experimentally
found that they are good estimates for the error of our numerical
methods.  Therefore, we use the relative errors in area and length to
study the convergence rates of our different time integrators.  We also
report minimum distances between vesicles and solid walls to demonstrate
convergence, but without an exact solution, we must used an over-refined
solution as the exact solution.  The experiments we perform are now
summarized.

\begin{itemize}
\item{\bf Shear Flow} 
(Tables~\ref{t:shear:RA65}--\ref{t:shear:RA99} and
Figure~\ref{f:shear:reducedAreas}): Here we consider the classic
problem of two vesicles in the shear flow $\uu = (y,0)$ with the left
vesicle slightly elevated from the right vesicle.  We study the effect
of the reduced area and the viscosity contrast on the dynamics.  Then,
we compare the three different stable time integrators we have
discussed: first-order with explicit and semi-implicit inter-vesicle
interactions and second-order with semi-implicit inter-vesicle
interactions.  All three methods converge at the expected rate.
Moreover, if we use first-order time stepping, semi-implicit and
explicit inter-vesicle interactions give similar errors.  Therefore,
for first-order time stepping, if the time step is small enough, using
explicit inter-vesicle interactions is appropriate since the
resulting linear system is easier to solve.

\item{\bf Taylor-Green Flow}
(Table~\ref{t:taylorGreen:convergence} and
Figures~\ref{f:taylorGreen:collision}--\ref{f:taylorGreen:summary}):
Here we look at the stability of the semi-implicit inter-vesicle
interactions.  For a short time horizon, we quantify its stability.  We
do a convergence study and also discuss a significant difference in the
behaviour of the first- and second-order time integrators.

\item {\bf Extensional Flow}
(Table~\ref{t:extensional:convergence} and
Figure~\ref{f:extensional:summary}): We consider two vesicles placed
symmetrically in an extensional flow.  We show that the simulation is
unstable without near-singular integration.  A convergence study shows
that we achieve the desired first- and second-order convergence.  We
also report CPU times to demonstrate that we are achieving the expected
complexity.

\item{\bf Stenosis} 
(Table~\ref{t:stenosis:errors}--\ref{t:stenosis:gapSize} and
Figure~\ref{f:stenosis:snaps}): We consider a single vesicle passing
through a constriction in a confined flow.  We are able to take a
vesicle that is over two times taller than a similar simulation
in~\cite{rah:vee:bir}.  A convergence study on the errors in area and
length as well as the minimium distance between the vesicle and solid
wall show that we achieve second-order accuracy.

\item{\bf Couette Apparatus}
(Figures~\ref{f:couetteLow}--\ref{f:couetteHigh}): High concentration
flows are demonstrated with two Couette simulations, one with a moderate
volume fraction and one with a high volume fraction.  We compare first-
and second-order time stepping with the moderate volume fraction and
observe a similar behaviour in the errors that is described in the
Taylor-Green example.  For the larger volume fraction, we report results
for first- and second-order time stepping which both require
semi-implicit inter-vesicle interactions for the chosen time step size
to be stable.

\end{itemize}

Other details of our results include:
\begin{itemize}
\item The parameters we investigate are the number of points per
vesicle, $N$, the number of vesicles, $M$, the number of points per
solid wall, $N_{\mathrm{wall}}$, the time step size, $\Delta t$, the
viscosity contrast, $\nu$, and the time integrator.

\item An important dimensionless quantity is the vesicle's reduced
area: the ratio of its area to the area of a circle of the same
perimeter.  Because of the incompressibility and inextensibility
constraints, the reduced area should remain constant.  Intuitively, the
reduced area indicates how much the vesicle can deform.  The closer the
reduced area is to its maximum value of one, the less it is able to
deform.  We are particularly interested in vesicles of reduced area
0.65. This value results in shapes that resemble (two-dimensional
approximations of) red blood cells.

\item All vesicle-to-vesicle and vesicle-to-wall interactions are
accelerated with the FMM.  The wall-to-wall iterations are precomputed
and stored and the wall-to-vesicle interactions are computed with the
direct method which results in a $\bigO(MNN_{\mathrm{wall}})$
calculation.  For flows with viscosity contrast, the double-layer
potential is also evaluated with the direct method resulting in a
$\bigO(M^{2}N^{2})$ calculation.  Fast algorithms for these
interactions are under consideration.

\item We precondition all the linear systems with the block-diagonal
preconditioner.  We have observed that this results in a
mesh-independent solver with semi-implicit interactions.

\item We have experimented with different GMRES tolerances.  If it is
too large, second-order convergence is difficult to achieve.  If it is
to small, the condition number of the block-diagonal preconditioner
makes the tolerance unachievable.  Experimentally we have found that
$10^{-12}$ is a good compromise.

\item The bending modulus sets the time scale.  We let $\kappa_{b}=0.1$
for all simulations.

\item All simulations are preformed in MATLAB on a six-core 2.67GHz
Intel Xeon processor with 24GB of memory. 
\end{itemize}

\paragraph{Estimation of the overall complexity}
Here we summarize the cost of the most expensive parts of the overall
algorithm.  Recall that we are using a mesh-independent preconditioned
GMRES iteration.  Therefore, the cost can be summarized by the time
required to preform one matrix-vector product.
\begin{itemize}
  \item{\em {Vesicle-vesicle interactions:}}
    Using the FMM, the single-layer potential requires $\bigO(MN)$
    operations. 

  \item{\em {Viscosity contrast:}}
    Since we do not have a fast algorithm to compute the double-layer
    potential, viscosity contrast requires $\bigO(M^{2}N^{2})$
    operations.  However, this is just a shortcoming of the specific
    implementation, not a theoretical one.

  \item{\em {Confined flows:}}
    Wall-to-wall interactions are precomputed and stored.
    Wall-to-vesicle interactions require $\bigO(MNN_{\mathrm{wall}})$
    operations and vesicle-to-wall interactions require
    $\bigO(MN+N_{\mathrm{wall}})$ operations using the FMM.

  \item{\em {Near-singular integration:}}
    Our near-singular integration complexity is summarized in
    Section~\ref{s:near-singular}.  The upsampling results in an
    increase in the above complexities.  For instance, vesicle-vesicle
    interactions with the FMM require $\bigO(MN^{3/2}+MN)$ operations
    and vesicle-to-wall interactions require
    $\bigO(MN^{3/2}+N_{wall})$.  However, since we have not implemented
    a fast summation for the double-layer potential, flows with a
    viscosity contrast require a $\bigO(M^{2}N^{5/2})$ calculation.
\end{itemize}

\subsection{Shear Flow}
Here we consider two vesicles in the shear flow $(y,0)$.  We initially
place the left vesicle's center slightly above the x-axis and place the
right vesicle's center exactly on the x-axis.  Initially, we expect the
right vesicle to undergo tank-treading~\cite{kantsler,misbah2006} while
the other vesicle will travel from left to right.  Here we study the
effect of the reduced area and the viscosity contrast on the dynamics
and errors.  We expect vesicles with larger viscosity contrast and
larger reduced areas to behave more like rigid bodies.  In particular,
we expect them to come closer to one another.
Figure~\ref{f:shear:reducedAreas} indicate that a higher viscosity
contrast results in the vesicles coming closer together.  Moreover, for
larger reduced areas, the effect of the viscosity contrast is less
pronounced.  We also see that vesicles with a larger reduced area come
closer together.

\begin{figure}[htps]
\centering
\begin{tabular}{m{1cm}CCCCC}
$\nu=1$ &
\input{shearOrder1RA1viscCont1Time1.tikz} &
\input{shearOrder1RA1viscCont1Time2.tikz} &
\input{shearOrder1RA1viscCont1Time3.tikz} &
\input{shearOrder1RA1viscCont1Time4.tikz} &
\input{shearOrder1RA1viscCont1Time5.tikz} \\
$\nu=4$ &
\input{shearOrder1RA1viscCont4Time1.tikz} &
\input{shearOrder1RA1viscCont4Time2.tikz} &
\input{shearOrder1RA1viscCont4Time3.tikz} &
\input{shearOrder1RA1viscCont4Time4.tikz} &
\input{shearOrder1RA1viscCont4Time5.tikz} \\
$\nu=1$ &
\input{shearOrder1RA2viscCont1Time1.tikz} &
\input{shearOrder1RA2viscCont1Time2.tikz} &
\input{shearOrder1RA2viscCont1Time3.tikz} &
\input{shearOrder1RA2viscCont1Time4.tikz} &
\input{shearOrder1RA2viscCont1Time5.tikz} \\
$\nu=4$ &
\input{shearOrder1RA2viscCont4Time1.tikz} &
\input{shearOrder1RA2viscCont4Time2.tikz} &
\input{shearOrder1RA2viscCont4Time3.tikz} &
\input{shearOrder1RA2viscCont4Time4.tikz} &
\input{shearOrder1RA2viscCont4Time5.tikz}
\end{tabular}
\begin{tabular}{cc}
\input{shearRA1GapSizeOrder1.tikz} &
\input{shearRA2GapSizeOrder1.tikz}
\end{tabular}
\mcaption{Two vesicles are submerged in a shear flow.  In the first and
second rows, the vesicles have reduced area 0.65, and in the third and
fourth rows, the vesicles have reduced area 0.99.  Each simulation is
done with viscosity contrasts $\nu=1$ and $\nu=4$.  The bottom plots
show the distance between the vesicles for both reduced areas and both
viscosity contrasts.  Comparing semi-implicit inter-vesicle interactions
with both first- and second-order time stepping, the distance between
the vesicles never differ by more than 2.2\% for both reduced areas and
both viscosity contrasts.}{f:shear:reducedAreas} 
\end{figure}

We now do a convergence study for both values of the reduced area and
viscosity contrast.  We vary the number of points per vesicle, the
number of time steps, and the time integrator.  We report the errors in
area and length at the time horizon $t=12$ in Tables~\ref{t:shear:RA65}
and~\ref{t:shear:RA99}.  First, we see that first- and second-order
convergence is achieved for all combinations of reduced area and
viscosity contrast.  Second, semi-implicit inter-vesicle interactions
do not improve the accuracy of first-order time stepping, but, we will
see later that this time integrator is more stable.

\begin{table}[htps]
\centering
\begin{tabular}{ccc|cc|cc} 
Time & & & \multicolumn{2}{c|}{$\nu=1$} & \multicolumn{2}{c}{$\nu=4$} \\
Integrator   & $N$ & $\Delta t$ & $e_{A}$   & $e_{L}$   & 
$e_{A}$   & $e_{L}$ \\ 
\hline

Explicit (1) & 32  & 0.04       & $6.73e-2$ & $5.21e-2$ &
$4.07e-2$ & $2.41e-2$ \\ 

Explicit (1) & 64  & 0.02       & $3.36e-2$ & $2.62e-2$ &
$2.02e-2$ & $1.20e-2$ \\ 

Explicit (1) & 128 & 0.01       & $1.68e-2$ & $1.31e-2$ &
$1.01e-2$ & $5.99e-2$ \\ 
\hline

Implicit (1) & 32  & 0.04       & $6.73e-2$ & $5.21e-2$ &
$4.06e-2$ & $2.40e-2$ \\ 

Implicit (1) & 64  & 0.02       & $3.36e-2$ & $2.62e-2$ &
$2.02e-2$ & $1.20e-2$ \\ 

Implicit (1) & 128 & 0.01       & $1.68e-2$ & $1.31e-2$ &
$1.02-2$ & $5.98e-3$ \\ 
\hline

Implicit (2) & 32  & 0.04       & $1.59e-4$ & $3.21e-4$ &
$2.69e-5$ & $6.33e-5$ \\ 

Implicit (2) & 64  & 0.02       & $3.87e-6$ & $6.52e-5$ &
$3.78e-6$ & $1.20e-5$ \\ 

Implicit (2) & 128 & 0.01       & $2.79e-7$ & $1.47e-5$ &
$3.85e-7$ & $2.28e-6$ 
\end{tabular}
\mcaption{The errors in area and length at $t=12$ for different time
integrators applied to two vesicles of reduced area $0.65$ in a shear
flow.  The parenthetic values are the time stepping order.  Results are
presented for no viscosity contrast and a viscosity contrast of
4.}{t:shear:RA65}

\begin{tabular}{ccc|cc|cc} 
Time & & & \multicolumn{2}{c|}{$\nu=1$} & \multicolumn{2}{c}{$\nu=4$} \\
Integrator   & $N$ & $\Delta t$ & $e_{A}$   & $e_{L}$   & 
$e_{A}$   & $e_{L}$ \\ 
\hline

Explicit (1) & 32  & 0.04       & $1.05e-1$ & $5.48e-2$ &
$9.67e-2$ & $4.94e-2$ \\ 

Explicit (1) & 64  & 0.02       & $5.15e-2$ & $2.70e-2$ &
$4.78e-2$ & $2.45e-2$ \\ 

Explicit (1) & 128 & 0.01       & $2.55e-2$ & $1.34e-2$ &
$2.37e-2$ & $1.22e-2$ \\ 
\hline

Implicit (1) & 32  & 0.04       & $1.05e-1$ & $5.48e-2$ &
$9.69e-2$ & $4.93e-2$ \\ 

Implicit (1) & 64  & 0.02       & $5.16e-2$ & $2.70e-2$ &
$4.77e-2$ & $2.45e-2$ \\ 

Implicit (1) & 128 & 0.01       & $2.55e-2$ & $1.34e-2$ &
$2.36e-2$ & $1.22e-2$ \\ 
\hline

Implicit (2) & 32  & 0.04       & $2.40e-5$ & $1.13e-4$ &
$5.85e-5$ & $4.68e-5$ \\ 

Implicit (2) & 64  & 0.02       & $1.65e-6$ & $1.50e-5$ &
$7.05e-6$ & $5.86e-6$ \\ 

Implicit (2) & 128 & 0.01       & $1.56e-7$ & $1.92e-6$ &
$8.62e-7$ & $7.26e-7$ 

\end{tabular}

\mcaption{The errors in area and length at $t=12$ for different time
integrators applied to two vesicles of reduced area $0.99$ in a shear
flow.  The parenthetic values are the time stepping order.  Results are
presented for no viscosity contrast and a viscosity contrast of
4.}{t:shear:RA99}
\end{table}

\subsection{Taylor-Green Flow}
\label{s:taylorGreen}
We consider the Taylor-Green flow $\uu = (\sin x \cos y,-\cos x \sin
y)$ with 9 vesicles each with $N=64$ points of reduced area 0.65 placed
inside the periodic cell $(0,2\pi)^{2}$ (left plot of
Figure~\ref{f:taylorGreen:collision}).  We consider first-order IMEX
Euler with both semi-implicit and explicit inter-vesicle interactions.
We find the largest time step of both methods so that the vesicles do
not cross before a time horizon of $T = 1$.  With explicit
interactions, the vesicles cross if $74$ time steps are taken but not
if $75$ time steps are taken.  However, with semi-implicit
interactions, the simulation successfully completes even with one time
step.  Of course the error is too large and a time step size this large
is unreasonable.

\begin{figure}[htps]
\begin{center}
  \begin{tabular}{ccc}
    \input{taylorGreenSnapsExplicitTime1.tikz} &
    \input{taylorGreenSnapsExplicitTime2.tikz} &
    \input{taylorGreenSnapsImplicitTime2.tikz} \\
    &
    \input{taylorGreenSnapsExplicitTime2Zoomed.tikz} &
    \input{taylorGreenSnapsImplicitTime2Zoomed.tikz}
  \end{tabular}
\end{center}
\mcaption{We compare first-order time stepping with explicit and
semi-implicit inter-vesicle interactions for a time step size of 0.02
and 64 points per vesicle with the initial configuration in the
left-most plot.  With explicit interactions (middle plot), the vesicles
cross after only 13 time steps.  However, with semi-implicit
interactions (right plot), the vesicles do not cross.  The bottom plots
show magnifications of the region where the explicit solver first
crosses.}{f:taylorGreen:collision}
\end{figure}

For this example, we found that using explicit interactions with 75
time steps requires approximately the same amount of CPU time as using
semi-implicit interactions with 30 time steps.  Therefore, if one
permits the error to be approximately 2.5 times larger, using
semi-implicit inter-vesicle interactions is justified.

We do a convergence study for first- and second-order time stepping
with semi-implicit inter-vesicle interactions.  If these interactions
are treated explicitly, the coarsest simulation results in the vesicles
crossing.  We achieve the desired first-order convergence, but
second-order convergence is not achieved at these resolutions.
However, we ran the simulation with $\Delta t = 3.13e-4$ and $\Delta t
= 1.56e-4$ and this is where second-order convergence is first
observed; the error in area reduced by $4.67$ and the error in length
reduced by $3.96$.

We have also observed a qualitative difference in the
behaviour of the errors of first- and second-order time stepping.  The
errors of the first-order methods continue to grow throughout the
simulation, while for second-order time stepping, the errors
plateau.  Thus, given a fixed time step size, second-order time stepping
can simulate longer time horizons than first-order time stepping.  We
demonstrate this behaviour in the bottom plot of
Figure~\ref{f:taylorGreen:summary} where we plot the errors in area and
length for the first- and second-order time stepping.  Snapshots of the
simulation are in the top plots of Figure~\ref{f:taylorGreen:summary}.

\begin{table}[htps]
\begin{center}
\begin{tabular}{cc|cc|cc}
& & \multicolumn{2}{c|}{First-order} & \multicolumn{2}{c}{Second-order} \\
$N$ & $\Delta t$ & $e_{A}$ & $e_{L}$ & $e_{A}$ & $e_{L}$ \\
\hline
32  & $2e-2$ & $5.24e-2$ & $3.03e-2$ & $5.25e-4$ & $3.98e-3$ \\
64  & $1e-2$ & $2.77e-2$ & $1.64e-2$ & $2.08e-4$ & $1.96e-3$ \\
128 & $5e-3$ & $1.43e-2$ & $8.68e-3$ & $7.76e-5$ & $9.42e-4$ \\
256 &$2.5e-3$& $7.24e-3$ & $4.53e-3$ & $2.73e-5$ & $4.20e-4$ \\
\end{tabular}
\mcaption{The errors in area and length at $t=5$ for a Taylor-Green
flow.  All runs were done with semi-implicit inter-vesicle interactions
since explicit interactions with $N=32$ and $\Delta t = 2e-2$ results
in crossing vesicles.}{t:taylorGreen:convergence} 
\end{center}
\end{table}

\begin{figure}[htps]
\begin{center}
  \begin{tabular}{cccccc}
    \input{taylorGreenSnapsTime1.tikz} &
    \input{taylorGreenSnapsTime2.tikz} &
    \input{taylorGreenSnapsTime3.tikz} &
    \input{taylorGreenSnapsTime4.tikz} &
    \input{taylorGreenSnapsTime5.tikz} &
    \input{taylorGreenSnapsTime6.tikz}
  \end{tabular}
  \begin{tabular}{cc}
    \input{taylorGreenErrorsOrder1.tikz} & 
    \input{taylorGreenErrorsOrder2.tikz} 
  \end{tabular}
\end{center}
\mcaption{The top plots show several snapshots of a Taylor-Green
simulation with a long time horizon using second-order time stepping.
In the bottom plots are the errors in area and length for the first- and
the second-order time integrators.  Notice that the errors for the
first-order method continue to grow and the simulation stops near $t=20$
(the simulation is stopped when a 10\% error has been incurred), whereas
the errors for the second-order method plateau and the simulation is
able to run much longer.}{f:taylorGreen:summary} 
\end{figure}

\subsection{Extensional Flow}
We consider two vesicles with reduced area $0.65$ symmetrically placed
around the origin (top plots of Figure~\ref{f:extensional:summary}).
The background velocity is given by $\uu = (-x,y)$.  We first run the
simulation with and without near-singular integration.  We compute the
distance between the vesicles, and the errors in area and length
(bottom plots of Figure~\ref{f:extensional:summary}).  We see that when
$t \approx 5$, the simulation without near-singular integration begins
to introduce large errors in the area of the vesicles.  This indicates
that the approximation of the single-layer potential is not adequate.
Then, around $t \approx 12$, the vesicles are too close for the
$N$-point trapezoid rule and the simulations fails.  This agrees with
our error estimate from Appendix~\ref{A:AppendixA}  since for this
simulation, an arclength term is approximately 0.24.  However, with
near-singular integration, the simulation successfully runs to
completion.

\begin{figure}[htps]
\begin{center}
  \begin{tabular}{ccccc}
  \input{extensionalSnapsTime1.tikz} &
  \input{extensionalSnapsTime2.tikz} &
  \input{extensionalSnapsTime3.tikz} &
  \input{extensionalSnapsTime4.tikz} &
  \input{extensionalSnapsTime5.tikz} 
  \end{tabular}
  \begin{tabular}{cc}
  \input{extensionalGapSize.tikz} &
  \input{extensionalErrors.tikz}
  \end{tabular}
\end{center}
\mcaption{Top: Two vesicles discretized with $N=32$ points at several
time steps using near-singular integration.  Bottom left: The distance
between the two vesicles both with (solid) and without (dashed)
near-singular integration.  Bottom right: The error in area (blue) and
length (red) with (solid) and without (dashed) near-singular
integration.}{f:extensional:summary} 
\end{figure}

We report the errors in area and length at $t=24$ in
Table~\ref{t:extensional:convergence}.  Also reported is the error in
the distance between the vesicles at $t=24$, for which we take the
``true'' distance from an overrefined simulation.  As before, we
achieve first-order convergence and there is little difference between
the accuracy when comparing explicit and semi-implicit inter-vesicle
interactions.  We also achieve second-order convergence, but it
requires a finer time step size when compared to the first-order
methods before the asymptotic rate is achieved.

In Table~\ref{t:extensional:convergence}, we also list the required CPU
time per time step relative to the smallest simulation.  When the
number of points is doubled from 32 to 64, the amount of work per time
step goes up by less than a factor of two.  The increase in CPU time
when $N$ is doubled from 64 to 128 is attributed to extra work required
by near-singular integration since, to simplify the implementation, we
only work with powers of 2.  Therefore, we upsample 64 to 256, but 128
is upsampled to 4096.  We also see that semi-implicit inter-vesicle
interactions are computationally not much more expensive.  This is
because a lot of work is spent precomputing the block-diagonal
preconditioner which is done once per time step for all the time
integrators.

\begin{table}[htps] 
\begin{centering} 
\begin{tabular}{ccccccc} 
Integrator & $N$ &   $\Delta t$ & $e_{A}$ & $e_{L}$ & $e_{\mathrm{gap}}$ & CPU Time \\ \hline
Explicit (1) & $32$  & $0.04$ & $1.99e-4$ & $6.53e-4$ & $1.97e-1$ & --- \\  
Explicit (1) & $64$  & $0.02$ & $6.34e-5$ & $3.42e-4$ & $3.73e-3$ & $1.15$ \\ 
Explicit (1) & $128$ & $0.01$ & $3.13e-5$ & $1.78e-4$ & $1.16e-3$ & $2.83$ \\ 
\hline
Implicit (1) & $32$  & $0.04$ & $2.12e-4$ & $6.64e-4$ & $1.92e-1$ & $2.65$ \\  
Implicit (1) & $64$  & $0.02$ & $6.43e-5$ & $3.46e-4$ & $2.23e-3$ & $2.72$ \\ 
Implicit (1) & $128$ & $0.01$ & $3.15e-5$ & $1.79e-4$ & $4.35e-4$ & $7.40$ \\ 
\hline
Implicit (2) & $32$  & $0.04$ & $3.26e-4$ & $1.38e-4$ & $3.96e-2$ & $2.43$ \\  
Implicit (2) & $64$  & $0.02$ & $4.14e-6$ & $4.45e-5$ & $8.32e-2$ & $2.89$ \\ 
Implicit (2) & $128$ & $0.01$ & $5.08e-7$ & $1.27e-5$ & $2.69e-2$ & $8.13$ \\
Implicit (2) & $128$ & $0.005$& $1.02e-7$ & $3.13e-6$ & $1.27e-4$ & $7.70$ \\
\end{tabular}
\mcaption{The errors in area, length, and the error in the gap size at
$t=24$, and the CPU time per time step for the extensional flow.  The
parenthetic values are the order of the time integrator.  The error in
the gap size is found by computing a ``true'' solution with $N=256$ points per
vesicle and a time step size so that the error in area is $9.0e-9$ and
in length is $7.3e-11$.  The CPU times are relative to the simulation
with $N=32$, $\Delta t = 0.04$ which took approximately
$2.62e-1$ seconds per time step.}{t:extensional:convergence}
\end{centering}
\end{table}

\subsection{Stenosis}

We consider a single vesicle of reduced area $0.65$ in a constricted
tube.  The solid wall is parameterized as $\xx(\theta) =
10r(\theta)\cos(\theta)$ and $y(\theta) =
3r(\theta)\sin(\theta)\eta(\theta)$, where $r(\theta) =
(\cos(\theta)^{8} + \sin(\theta)^{8})^{-1/8}$, and
\begin{align*}
  \eta(\theta) = \left\{
    \begin{array}{cl}
      \frac{1-0.6\cos(x)}{1.6} & |x(\theta)| \leq \pi, \\
      1 & |x(\theta)| > \pi.
    \end{array}
  \right.
\end{align*}
For this geometry, the error of the trapezoid rule is $\bigO(10^{-4})$
with $N_{\mathrm{wall}}=256$ points and is $\bigO(10^{-6})$ with
$N_{\mathrm{wall}}=512$ points.  We impose a Gaussian-shaped boundary
condition at the intake and outake of the tube which is plotted along
with the geometry in Figure~\ref{f:stenosis:snaps}.  The vesicle's
initial height is $2.3$ times larger than the size of the constriction.
We report a convergence study for second-order time stepping in
Table~\ref{t:stenosis:errors} and we exceed the expected rate of
convergence.  We note that for this example, $N_{\mathrm{wall}}$ can not
be less than 256.  If it is, the Lagrange interpolation points required
by near-singular integration for wall-to-vesicle interactions are on
both sides of the solid wall.  Since the double-layer potential has a
jump across the solid wall, this creates a discontinuity in the velocity
field and the method becomes unstable.

This example is ideally situated for adaptive time stepping.  As the
vesicle passes through the constriction, large jumps in the errors are
committed.  In the simulation with $\Delta t = 5.0e-3$, the error in
length jumps by four orders of magnitude, and the error in area jumps
by three orders of magnitude near $t=3$.  These jumps are less severe
in the simulation with $\Delta t = 2.5e-3$, which has a jump of one
order of magnitude in length and two orders of magnitude in area.

\begin{figure}[htps] 
\centering 
\begin{tabular}{ccc}
\input{stenosisSnapsTime1.tikz} &
\input{stenosisSnapsTime2.tikz} &
\input{stenosisSnapsTime3.tikz} \\
\input{stenosisSnapsTime4.tikz} &
\input{stenosisSnapsTime5.tikz} &
\input{stenosisSnapsTime6.tikz} \\
\input{stenosisSnapsTime7.tikz} &
\input{stenosisSnapsTime8.tikz} &
\input{stenosisSnapsTime9.tikz}
\end{tabular} 
\mcaption{Snapshots of a single vesicle in a constricted tube.  The
boundary condition on the solid wall (plotted in blue) at the inlet and
outlet is Gaussian-shaped, and on the remainder of the solid wall is
zero.  The vesicle is initially $2.3$ times larger than the size of the
constriction.  Results from~\cite{rah:vee:bir} were only able to
simulate vesicles that were $1.12$ times larger than the
constriction.}{f:stenosis:snaps}
\end{figure}

\begin{table}[htps] 
\begin{centering} 
  \begin{tabular}{cccccc} 
    $N$ & $N_{\mathrm{wall}}$ & $\Delta t$ & 
      $e_{A}$ & $e_{L}$ \\ \hline 
    128 & 256 & $5.0e-3$ & $6.55e-4$ & $1.17e-2$ \\ 
    256 & 512 & $2.5e-3$ & $9.59e-7$ & $1.05e-6$  
  \end{tabular} 
  \mcaption{The errors in area and length at $t=8$ for the second-order
  stenosis flow.}{t:stenosis:errors}
  \end{centering}
\end{table}

We also report the minimum distance, and the time when the minimum
distance occurs, between the vesicle and the solid wall in
Table~\ref{t:stenosis:gapSize}.  In order to focus the error on the time
step size, we take a fine spatial resolution of $N=256$ and
$N_{\mathrm{wall}}=512$.  We report the minimum over time of the
distance 
\begin{align*}
  d(\gamma,\Gamma) = \inf_{\xx \in \gamma, \yy \in \Gamma} \|\xx - \yy\|
\end{align*}
rather than the error since an exact solution is not available.  The
time integrators predict different locations, $t_{\min}$, where the
minimum occurs, but they both are converging.  For the first-order
results, if we compare the minimum distances with the distance of the
most accurate second-order simulation at $t=5.26$, which is $9.6369e-2$,
then we see that first-order convergence has been achieved.  Without a
much more accurate numerical solution, we are unable to compute
convergence rates for the second-order results.  However, it appears
that we have resolved the distance up to at least three digits of
accuracy.

\begin{table}[htps] 
  \begin{centering} 
  \begin{tabular}{c|cc|cc} 
  & \multicolumn{2}{c|}{First-order} &
    \multicolumn{2}{c}{Second-order} \\
    $\Delta t$ & $t_{\min}$ & $d(\Gamma,\gamma)$ &  
      $t_{\min}$ & $d(\Gamma,\gamma)$  \\ \hline
    $1.0e-2$ & $4.37$ & $1.0318e-1$ & $4.34$ & $1.0196e-1$ \\
    $5.0e-3$ & $4.41$ & $1.0262e-1$ & $4.39$ & $9.8055e-2$ \\
    $2.5e-3$ & $4.42$ & $1.0143e-1$ & $5.77$ & $9.5809e-2$ \\
    $1.3e-3$ & $5.26$ & $9.9541e-2$ & $5.77$ & $9.5804e-2$ \\
    $6.3e-4$ & $5.26$ & $9.7938e-2$ & $5.77$ & $9.5803e-2$ \\ 
    $3.1e-4$ & $5.26$ & $9.7137e-2$ & $5.77$ & $9.5802e-2$ 
  \end{tabular}
 \mcaption{The minimum distance, $d(\Gamma,\gamma)$, between the vesicle
 and the solid wall, and the time, $t_{\min}$, when this distance is
 achieved.  Rather than forming an overrefined ''true" solution, we
 report the actual distances.}{t:stenosis:gapSize}
 \end{centering}
\end{table}

\subsection{Couette Apparatus}
\label{s:couette}
Here we consider two simulations of a Couette apparatus with different
volume fractions.  We first consider a Couette apparatus where the
inner boundary, which is slightly off-centered, is rotating with
constant angular velocity and the outer boundary is fixed.  We randomly
place 42 vesicles of reduced area 0.65 inside the apparatus which
corresponds to a volume fraction of 20\%.  Each vesicle is discretized
with 64 points, each solid wall with 64 points, the time step size is
$0.01$, and the time horizon is $T=30$.  Snapshots of the simulation
are in the top plots of Figure~\ref{f:couetteLow}.  Using first-order
time stepping with explicit inter-vesicle interactions, the error in
area is $1.70e-2$ and the error in length is $7.94e-2$ at $t=30$.
Semi-implicit inter-vesicle interactions are not required for this time
step size.  The simulation took a little less than 6 hours, or about
7.1 seconds per time step.

We repeat this experiment with second-order time stepping which
requires semi-implicit inter-vesicle interactions.  The simulation took
about 83 hours, or about 100 seconds per time step.  The errors in area
and length at $t=30$ are $2.95e-3$ and $7.13e-4$.  To achieve a similar
accuracy with first-order time stepping, a time step size 11 times
smaller would be required which would require more than 83 hours of CPU
time due to the double-layer potentials.  Moreover, the errors in area
and length for second-order time stepping have plateaued while for
first-order time stepping, they continue to grow (bottom plots of
Figure~\ref{f:couetteLow}).  Again, this means that with second-order
time stepping, we can take much longer time horizons without requiring
a smaller time step size.

\begin{figure}[htps]
\begin{center}
  \begin{tabular}{c@{\,}c@{\,}c@{\,}c@{\,}c@{\,}c@{\,}}
    \input{couetteSnapsTime1.tikz} &
    \input{couetteSnapsTime2.tikz} &
    \input{couetteSnapsTime3.tikz} &
    \input{couetteSnapsTime4.tikz} &
    \input{couetteSnapsTime5.tikz} &
    \input{couetteSnapsTime6.tikz} 
  \end{tabular}
  \begin{tabular}{cc}
    \input{couetteErrorsOrder1.tikz} & 
    \input{couetteErrorsOrder2.tikz} 
  \end{tabular}
\end{center}
\mcaption{The top plots show snapshots of a Couette flow with 42
vesicles and a volume fraction of 20\%.  One vesicle is coloured in blue
to view its time history.  The outer boundary is fixed while the inner
boundary completes one revolution every 10 time units.  The bottom plots
show the errors in area and length for the first- and second-order time
integrators.  Notice that the errors for first-order method  continue to
grow whereas the errors for the second-order method
plateau.}{f:couetteLow}
\end{figure}

We now consider a Couette apparatus with the inner boundary exactly in
the middle of the outer boundary.  We use 150 vesicles of reduced area
0.65 resulting in a volume fraction of approximately $40.5\%$.  With
this higher concentration, a finer discretization of the solid walls is
required to resolve its density function $\eeta$.  Therefore, we have
discretized both boundaries with 256 points.  We also used our
collision detection algorithm to alert us if vesicles had crossed.
With $\Delta t = 0.01$ and explicit inter-vesicle interactions, the
first-order method results in crossing vesicles.  Therefore, we report
the results using first-order semi-implicit inter-vesicle interactions
which does not result in vesicles crossing.  Snapshots are illustrated
in Figure~\ref{f:couetteHigh}.  On average, each time step required 27
GMRES iterations. The total simulation took about 10 days, but this
number will drop dramatically once the double-layer potential is
implemented with a fast summation method.  At $t=10$, the error in area
is $9.62e-3$ and the error in length is $1.92e-2$.  We also ran this
simulation with second-order time stepping and achieved an error in
area of $9.39e-4$ and an error in length of $4.19e-4$.
\begin{figure}[htps]
\begin{center}
  \begin{tabular}{c@{\,}c@{\,}c@{\,}c@{\,}c@{\,}c@{\,}}
    \input{couetteDenseSnapsTime1.tikz} &
    \input{couetteDenseSnapsTime2.tikz} &
    \input{couetteDenseSnapsTime3.tikz} &
    \input{couetteDenseSnapsTime4.tikz} &
    \input{couetteDenseSnapsTime5.tikz} &
    \input{couetteDenseSnapsTime6.tikz} 
  \end{tabular}
\end{center}
\mcaption{Snapshots of a Couette flow with 150 vesicles.  One vesicle
is coloured in blue to view its time history.  The inner boundary
completes one full revolution every 10 time units.}{f:couetteHigh}
\end{figure}

%% file: conclusions.tex
We have presented and tested a collection of extensions of the boundary
integral equation formulation for vesicle suspensions outlined
in~\cite{rah:vee:bir}.  Our goal is to create a robust method for
handling high concentration suspensions.  The main contributions we have
presented are:
\begin{itemize}

\item To remove stiffness, we have introduced a new time integrator that
treats inter-vesicle interactions semi-implicitly.  This allows us to
take larger time steps and to use second-order time stepping.  

\item To handle vesicles in close proximity, we employed and tested the
near-singular integration algorithm outlined
in~\cite{ying-biros-zorin06}.  This algorithm creates a uniform error
when evaluating integral operators, is easy to implement, and does not
significantly increase the overall complexity.  The vesicle flow
configurations that we tested are unstable without near-singular
integration.

\item To test for collisions, we introduced a spectrally accurate
algorithm that uses standard potential theory and is compatible with
our near-singular integration scheme and the fast multipole method.

\item We computed the pressure and stress of the single- and
double-layer potentials.  To use near-singular integration, we require
the limiting values of the pressure and stress as the target point
approaches the boundary.  All these jumps were computed.

\item We obtained first- and second-order convergence for a variety of
bounded and unbounded flows.  We also observed examples where the errors
for second-order convergence plateau while the error for first-order
convergence continues to grow.
\end{itemize}

While these contributions are a major step towards creating a robust
solver for high concentration vesicle suspensions, there are several key
features that are necessary and are currently being investigated.

\begin{itemize}

\item The vesicle shapes we have presented in Section~\ref{s:results}
can be well represented with $N=128$ or fewer points per vesicle.
However, we may require more points to represent quantities such as the
traction jump or the velocity field.  Therefore, spatial adaptivity is 
under consideration.

\item Currently, the time step is found using a trial and error
process.  We are developing high-order time integrators that use error
control to automatically and adaptively adjust the time step size.

\item Our formulation is two-dimensional.  The algorithms we have
introduced naturally extend to three-dimensions, but the linear systems
are far too expensive to solve without suitable preconditioners. 

\item Our formulation uses the steady Stokes equations for both the
fluid in the bulk and in the vesicle interior. If transient, inertial,
or viscoelastic effects are important one has to use a different
formulation.

\end{itemize}


%% file: appen1.tex
We compute error bounds for the near-singular integration scheme to
justify our choice of the $N^{3/2}$-point trapezoid rule for points in
the far zone.  Recall that the far zone is defined as all points that
are more than distance $h$ from a vesicle.  Consider a vesicle with
boundary $\gamma \in C^{\infty}$ and a density function $\ff \in
C^{M}$.  Recall that the single-layer potential is
\begin{align*}
  \uu(\xx) = \frac{1}{4\pi} \int_{\gamma} \left(-\log \rho + 
    \frac{\rr \otimes \rr}{\rho^{2}}\right)\ff(\yy) ds_{\yy}.
\end{align*}
Following~\cite{ying-biros-zorin06}, we need to bound the derivatives
of the kernel.  Bounding the derivatives of the logarithm term is
immediate, and the derivatives of the second term are bounded by noting
that
\begin{align*}
  \frac{\rr \otimes \rr}{\rho^{2}} = \frac{1}{\rho^{2}}S,
\end{align*}
where $S \in C^{\infty}$.  Bounds for the derivatives are
\begin{align*}
  |\partial^{\beta}_{\yy} \log \rho| \leq \frac{C}{\rho^{|\beta|+1}},   
  \hspace{20pt}
  \left|\partial^{\beta}_{\yy} \frac{\rr \otimes \rr}{\rho^{2}}\right| \leq
  \frac{C}{\rho^{|\beta|+2}},
\end{align*}
where $\beta$ is a multi-index of order $|\beta| \leq M$.  If $\rho >
\sqrt{h}$, then the $M^{th}$-order derivatives are bounded by
$Ch^{-M/2-1}$, and the error introduced by the $N$-point trapezoid rule
satisfies
\begin{align*}
  \epsilon_{trap} \leq h^{M} \cdot C h^{-M/2-1} = Ch^{M/2-1}.
\end{align*}
We would like to achieve the save accuracy when $\rho \in
(h,\sqrt{h}]$.  If $\rho \in (h,\sqrt{h}]$, the $M^{th}$-order
derivatives are bounded by $C h^{-M-2}$, and the error introduced by
the $N^{3/2}$-point trapezoid rule satisfies
\begin{align*}
  \epsilon_{trap} \leq h^{3M/2} \cdot C h^{-M-2} = Ch^{M/2-2}.
\end{align*}
Therefore, in what we have defined as the far zone, the error in the
trapezoid rule is $\bigO(h^{M/2-2})$.

To bound the error in the far zone of the double-layer potential, we note that
\begin{align*}
  \left|\partial^{\beta}_{\yy} \frac{\rr \cdot \nn}{\rho^{2}}
    \frac{\rr \otimes \rr}{\rho^{2}}\right| \leq \frac{C}{\rho^{|\beta|+4}}.
\end{align*}
We repeat the same calculation done for the single-layer potential and
the result is that the error in the far zone of the $N^{3/2}$-point
trapezoid rule applied to the double-layer potential is
$\bigO(h^{M/2-4})$.

For the near zone, we find the solution using $m+1$ Lagrange
interpolation points.  All the interpolation points except for the
first lie in the far zone and the above error bounds hold at these
points.  The Lagrange interpolation point on $\gamma$ is computed
using $N_{\mathrm{int}}$ interpolation points, so its error is
$\bigO\left(h^{N_{\mathrm{int}}-1}\right)$.  Combining the three
approximations, the error of our near-singular integration scheme for
the single-layer potential is
\begin{align*}
  \bigO\left(h^{\min\left(N_{\mathrm{int}}-1,m,\frac{M}{2}-2\right)}\right),
\end{align*}
and for the double-layer potential is
\begin{align*}
  \bigO\left(h^{\min\left(N_{\mathrm{int}}-1,m,\frac{M}{2}-4\right)}\right).
\end{align*}

%% file: appen2.tex
We compute the jumps in the pressure and the stress of the single-layer potential.  Consider a vesicle $\omega$ with boundary $\gamma$, normal vector $\nn$, tangent vector $\ttau$, and traction jump $\ff$.  Then, the pressure due to the single-layer potential is 
\begin{align*}
  p(\xx) &= \frac{1}{2\pi}\int_{\gamma} 
    \frac{\rr \cdot \ff}{\rho^{2}}ds_{\yy} \\
  &= \frac{1}{2\pi}\int_{\gamma} 
    \frac{(\rr \cdot \nn)}{\rho^{2}}(\ff \cdot \nn) ds_{\yy} + 
    \frac{1}{2\pi}\int_{\gamma} 
    \frac{(\rr \cdot \ttau)}{\rho^{2}}(\ff \cdot \ttau) ds_{\yy}.
\end{align*}
Both integrals are layer potentials of Laplace's equation.  The first is
the double-layer potential with density $\ff \cdot \nn$ and the second
is the adjoint of the tangential derivative of the single-layer
potential with density $\ff \cdot \ttau$.  Standard potential
theory~\cite{kellogg} shows that the first integral has a jump of size
$1/2$ while the second integral has no jump.  Therefore, 
\begin{align*}
  &\lim_{\substack{\xx \rightarrow \xx_{0} \\ \xx \in \omega}}p(\xx) = 
  +\frac{1}{2}\ff_{0} \cdot \nn_{0} + p(\xx_{0}), \\
  &\lim_{\substack{\xx \rightarrow \xx_{0} \\ \xx \notin \omega}}p(\xx) = 
  -\frac{1}{2}\ff_{0} \cdot \nn_{0} + p(\xx_{0}),
\end{align*}
where $\ff_{0} = \ff(\xx_{0})$ and $\nn_{0} = \nn(\xx_{0})$.

To find the jump in the stress tensor, we decompose it into its tangential and normal components
\begin{align*}
  T[\ssigma](\xx) &= \frac{1}{\pi}\int_{\gamma}\frac{\rr \otimes \rr}{\rho^{2}}
    \frac{\rr \cdot \ssigma}{\rho^{2}} \ff ds \\
  &= \frac{1}{\pi}\int_{\gamma}\frac{\rr \otimes \rr}{\rho^{2}}
    \frac{\rr \cdot \nn}{\rho^{2}}(\ssigma \cdot \nn)\ff ds + 
    \frac{1}{\pi}\int_{\gamma}\frac{\rr \otimes \rr}{\rho^{2}}
    \frac{\rr \cdot \ttau}{\rho^{2}}(\ssigma \cdot \ttau)\ff ds \\
    &=: T_{\nn}[(\ssigma \cdot \nn)\ff](\xx) + 
    T_{\ttau}[(\ssigma \cdot \ttau)\ff](\xx). 
\end{align*}
We use the usual method of placing a target point $\xx_{0} \in \gamma$, deforming the boundary by adding a circle of radius $\epsilon$ at $\xx_{0}$, computing the integral around the circle, and letting $\epsilon$ tend to zero.  However, we first rotate the tensors $T_{\nn}$ and $T_{\ttau}$ by the matrix $R$ so that the normal and tangential vectors at $\xx_{0}$ are $(1,0)$ and $(0,1)$ respectively.  Using the Residue Theorem to compute the contribution from the half circle, the jumps are 
\begin{align}
  &\lim_{\substack{\xx \rightarrow \xx_{0} \\ \xx \notin \omega}} 
    R T_{\nn}[\ff](\xx) R^{T} = \frac{1}{2}\left[
    \begin{array}{cc}
      1 & 0 \\ 0 & 1
    \end{array}
    \right] \ff(\xx_{0}) + R T_{\nn}[\ff](\xx_{0}) R^{T}, 
    \label{e:stress:norm:jump} \\
  &\lim_{\substack{\xx \rightarrow \xx_{0} \\ \xx \notin \omega}} 
    R T_{\ttau}[\ff](\xx) R^{T} = \frac{1}{2}\left[
    \begin{array}{cc}
      0 & 1 \\ 1 & 0
    \end{array}
    \right] \ff(\xx_{0}) + R T_{\ttau}[\ff](\xx_{0}) R^{T}. 
    \label{e:stress:tang:jump}
\end{align}
If the limits are taken with $\xx \in \omega$, the sign in the jump
changes to $-1/2$.  By left multiplying~\eqref{e:stress:norm:jump}
and~\eqref{e:stress:tang:jump} by $R^{T}$ and right multiplying by $R$,
the jump in $T[\sigma](\xx)$ from the interior of the vesicle
is~\eqref{e:stress:jump:int} and from the exterior of the vesicle
is~\eqref{e:stress:jump:ext}.

%% file: appen3.tex
The traction jump $\ff$ minimizes the energy functional
$\int_{\gamma}\frac{\kappa_{b}}{2}\kappa^{2}ds$ subject to the
inextensibility constraint.  An alternative traction jump is derived by
minimizing
\begin{align}
  \int_{\gamma}\frac{\kappa_{b}}{2}(\kappa - \tilde{\kappa})^{2}ds, \quad \text{subject to} \quad \|\xx_{s}\|=1,
    \label{e:new:bending}
\end{align}
where $\tilde{\kappa}$ is the curvature of a prescribed intrinsic
vesicle shape whose perimeter we are assuming is one.  By taking the
first variation of~\eqref{e:new:bending}, the force due to bending is
\begin{align*}
  \mathbf{f}_{b} = \left(\kappa_{ss} + \frac{\kappa^{3}}{2} -
    \tilde{\kappa}_{ss} - \frac{\kappa\tilde{\kappa}^{2}}{2}\right)\mathbf{n} + 
    (\kappa-\tilde{\kappa})\tilde{\kappa}_{s}\mathbf{t}.
\end{align*}
It is shown in~\cite{shravan} that terms of the form
$(h(s)\xx_{s})_{s}$ do no work for inextensible vesicles and hence
such forces can be added to modify $\mathbf{f}_{b}$.  A convenient
choice is $h(s) = 3/2\, \kappa^{2} - 1/2\, \tilde{\kappa}^{2}$ which
results in the simpler traction jump
\begin{align*}
  \ff = \ff_{b} + \ff_{\sigma} = \kappa_{b}(-\xx_{ssss} - \tilde{\kappa}_{ss} \mathbf{n} - 
    \kappa \tilde{\kappa}_{s}\mathbf{t}) + (\sigma \xx_{s})_{s}.
\end{align*}

We present simulations for several choices of $\tilde{\kappa}$ in
Figure~\ref{f:newbending}.  The flow is only driven by the vesicle
configuration since there is no background flow.  For $\tilde{\kappa}$,
we take the curvature of an ellipse with an aspect ratio of 1 (top
plots), 2 (middle plots), and 4 (bottom plots).  We see that the
curvature of the vesicle tries to fit to the curvature of the intrinsic
shape while maintaining its area and arclength.

\begin{figure}[htps]
\begin{center}
  \begin{tabular}{c@{\,}c@{\,}c@{\,}c@{\,}c@{\,}c@{\,}}
    \input{newBendingAspect1Time1.tikz} &
    \input{newBendingAspect1Time2.tikz} &
    \input{newBendingAspect1Time3.tikz} &
    \input{newBendingAspect1Time4.tikz} &
    \input{newBendingAspect1Time5.tikz} &
    \input{newBendingAspect1Curvature.tikz} \\

    \input{newBendingAspect2Time1.tikz} &
    \input{newBendingAspect2Time2.tikz} &
    \input{newBendingAspect2Time3.tikz} &
    \input{newBendingAspect2Time4.tikz} &
    \input{newBendingAspect2Time5.tikz} &
    \input{newBendingAspect2Curvature.tikz} \\

    \input{newBendingAspect4Time1.tikz} &
    \input{newBendingAspect4Time2.tikz} &
    \input{newBendingAspect4Time3.tikz} &
    \input{newBendingAspect4Time4.tikz} &
    \input{newBendingAspect4Time5.tikz} &
    \input{newBendingAspect4Curvature.tikz}
  \end{tabular}
\end{center}
\mcaption{Here we show the relaxation dynamics of a single vesicle with
three different prescribed intrinsic curvatures $\tilde{\kappa}$. There
is no imposed background flow.  The vesicle is discretized with $N=256$
points and $\Delta t = 0.02$.  After 200 time steps, the vesicle is
within $0.01 \%$ of its steady state configuration.  The final
curvature (red) and curvature of the intrinsic shape (blue) are in the
right plots.}{f:newbending} 
\end{figure}